\documentclass[12pt]{article}
\usepackage[english]{babel}
\usepackage[cp1250]{inputenc}
\usepackage{amsfonts, amsmath, amssymb}
\usepackage{amsthm}
\usepackage[pdftex]{geometry,graphicx,color}
\usepackage{tikz-cd}

\newtheorem{thm}{Theorem}
\newtheorem{prop}[thm]{Proposition}

\newtheorem{cor}[thm]{Corollary}
\newtheorem{rmk}[thm]{Remark}

\theoremstyle{definition}

\DeclareMathOperator{\Sym}{Sym}

\DeclareMathOperator{\x}{\tt x}
\DeclareMathOperator{\y}{\tt y}

\DeclareMathOperator{\PGL}{PGL}
\DeclareMathOperator{\Ia}{\mathbf{I}}
\DeclareMathOperator{\Ib}{\mathbf{II}}
\DeclareMathOperator{\Ic}{\mathbf{III}}

\begin{document}

\title{New examples of extremal positive linear maps}

\author{
Anita Buckley\thanks{%
Faculty of Mathematics and Physics, University of Ljubljana, Jadranska 19, SI-1000 Ljubljana, Slovenia, {\tt anita.buckley@fmf.uni-lj.si}. This author
was partially supported by the Slovenian Research Agency (J1-8132).}
\and
Klemen \v Sivic \thanks{%
Faculty of Mathematics and Physics, University of Ljubljana, Jadranska 19, SI-1000 Ljubljana, Slovenia, {\tt klemen.sivic@fmf.uni-lj.si}. This author
was partially supported by the Slovenian Research Agency (J1-8132, N1-0061, P1-0222).}
}

\date{}

\maketitle

\begin{abstract}
New families of  nonnegative biquadratic forms that have 8, 9 or 10 real zeros in $\mathbb{P}^2\times \mathbb{P}^2$ are constructed. These are the first examples with 8, 9 or 10 real zeros. It is known that nonnegative biquadratic forms with finitely many real zeros can have at most 10 zeros; our examples show that the upper bound is obtained. Such biquadratic forms define positive linear maps on real symmetric $3\times 3$ matrices that are not completely positive. Our constructions are explicit, and moreover we are able to determine which of the examples are extremal.  We extend the examples to positive maps on complex matrices and find families of extreme rays in the cone of positive maps.\\


{\bf Math. subj. class.:} 13J30, 15A86, 15B48, 47L07, 52A40
\end{abstract}

\section{Introduction}

\textbf{Notation}\hspace{3mm}
A linear map $\Phi \colon M_n(\mathbb{C})\to M_m(\mathbb{C})$ is \textit{positive} if it maps positive semidefinite matrices to positive semidefinite matrices, and is \textit{completely positive} (cp) if each of the ampliations 
$$\begin{array}{clcl} \vspace{2mm}
\Phi^{(k)}=I_k\otimes \Phi \colon &M_{kn}(\mathbb{C})=M_k(\mathbb{C})\otimes M_n(\mathbb{C})&
\to& M_{km}(\mathbb{C}) \\ 
 & [X_{ij}]_{i,j=1}^k&\mapsto &  [\Phi (X_{ij})]_{i,j=1}^k
\end{array}$$
is positive. 

We approach the study of positive linear maps from the viewpoint of real algebraic geometry. Therefore we restrict ourselves to positive linear maps between real symmetric matrices $\Phi \colon \Sym_n\to \Sym_m$. (Throughout the paper we will only write the upper triangle of a symmetric matrix). Such map can be easily extended to a positive map  
$\Phi \colon M_n(\mathbb{C})\to M_m(\mathbb{C})$ as shown in Section~\ref{SecExtension} and in \cite[Proposition 4.1]{KMcCSZ}.

The space of linear maps $\Phi \colon \Sym_n\to \Sym_m$ is isomorphic to the space of biquadratic forms in $n+m$ variables via the isomorphism 
$$\begin{array}{ccl} \vspace{2mm}
\Phi \mapsto p_{\Phi}(\x,\y)=\y^T\Phi \left( \x \x^T \right) \y, & \mbox{ where } &    
 \x = (x_0,\ldots,x_{n-1})\in \mathbb{P}^{n-1}(\mathbb{R}) \\
 &  \mbox{ and }  &  \y =(y_0,\ldots, y_{m-1}) \in \mathbb{P}^{m-1}(\mathbb{R}).
 \end{array}$$
As shown by Choi~\cite{C}, $\Phi$ is positive if and only if the polynomial $p_{\Phi}$ is nonnegative on $ \mathbb{P}^{n-1}(\mathbb{R})\times  \mathbb{P}^{m-1}(\mathbb{R})$, and $\Phi$ is completely positive if and only if $p_{\Phi}$ is a sum of squares (SOS) of bilinear forms. 

There is also a natural equivalence between positive linear maps on real symmetric matrices  and  quadratic matrix polynomials that are everywhere positive semidefinite (in other words,  semidefinite quadratic determinantal representations of  nonnegative polynomials). The equivalence is obtained by evaluating $\Phi$ on rank 1 symmetric matrices. It is easy to verify that 
$$\Phi \colon \!\! \Sym_n \! \! \to \! \Sym_m\,  \mbox{is positive} \Longleftrightarrow
\Phi \left(\x \x^T \right) \mbox{ is positive semidefinite for all}  \x \in \mathbb{P}^{n-1}(\mathbb{R}),$$
and $\Phi$ is completely positive if and only if 
$$\Phi \left(\x \x^T \right)=L(\x)L(\x)^T =\sum_{i=1}^l B_i \x \x^T\! B_i^T$$ 
for some $l\in\mathbb{N},$ some $m\times l$ linear matrix $L(\x)$  and some constant  $m\times n$ matrices $B_i$.
Completely positive map $\Phi$ is a \textit{congruence map} when $l=1$.

A positive map on real symmetric matrices (resp. nonnegative biquadratic form) is \textit{extremal} if it can not be written as a sum of two linearly independent positive maps (resp. nonnegative biquadratic forms). The definition of  \textit{extremal positive maps} from $M_n(\mathbb{C})$ to $M_m(\mathbb{C})$ is analogous. The convex cone of positive maps is the convex hull of extremal positive maps. \\

\noindent
\textbf{Background}\hspace{3mm}
This paper considers examples where $m=n=3$. The first example of a positive map on real symmetric matrices that is not completely positive is due to Choi~\cite{C},
$$\Phi\left(
\left[\!\! \begin{array}{ccc}
 z_{00}  & z_{01} & z_{02} \\
&  z_{11} &   z_{12} \\
 &  & z_{22} 
  \end{array}\!\! \right]
  \right)=
\left[\!\! \begin{array}{ccc}
 z_{00} +   2 z_{11} & -z_{01} & - z_{02} \\
&  z_{11} + 2 z_{22} &   -z_{12} \\
 &  & 2 z_{00} + z_{22} 
  \end{array}\!\! \right].$$
It is known that there are many more positive than completely positive maps~\cite{E, KMcCSZ, SWZ}. However, the examples in ~\cite{KMcCSZ} are not extremal. An older construction of Terpstra~\cite{ter} yields only special isolated examples which are in most cases not extremal. Most of the examples that can be found in the literature are obtained as generalizations of Choi's example, see \cite{CKL,  CK, H, Hou,  LF, O, Q, SC, TT}. 

By analogy with Blekherman's geometric explanation~\cite{Bleck} of the containment of the convex cone $\Sigma_{3,6}$ of SOS polynomials in the convex cone $P_{3,6} \subset \mathbb{P}^{27}(\mathbb{R})$ of nonnegative ternary polynomials of degree 6, we can study the containment of  the convex cone $\Sigma_{B(3,3)}$ of SOS biquadratic forms  in the convex cone  $P_{B(3,3)}$ of nonnegative biquadratic forms on $\mathbb{P}^2 \times \mathbb{P}^2$, or equivalently the containment of the convex cone of completely positive maps in the convex cone of positive maps from $\Sym_3$ to $\Sym_3$.
Determinant is a natural map 
$$\begin{array}{cccc}
\det: &  P_{B(3,3)} & \rightarrow &  P_{3,6} \\
        & p_{\Phi}(\x,\y) & \mapsto & \det \Phi(\x \x^T)
\end{array}$$  
that is obviously not injective. We also strongly believe that $\det$ is not surjective, however there exists no proof that some nonnegative polynomial, for example Robinson polynomial 
\begin{equation} \label{eq:therob} x_0^6+x_1^6+x_2^6
-x_0^4x_1^2-x_1^4x_2^2-x_2^4x_0^2- x_0^2x_1^4-x_1^2x_2^4-x_2^2x_0^4+3 x_0^2x_1^2x_2^2, \end{equation}
has no semidefinite quadratic determinantal representation.
It is not hard to see that $\det$ maps $\Sigma_{B(3,3)} $ into $\Sigma_{3,6} $, since for a completely positive map $\Phi$, the determinant $ \det \Phi(\x \x^T)$  is an SOS polynomial~\cite{Q}. On the other hand, Quarez' example~\cite[Proposition 5.1]{Q} with determinant 
$$\det \left[\!\! \begin{array}{ccc}
 x_0^2 +   x_2^2 & 0 & - x_0 x_2 \\
&  x_0^2 +   x_1^2 &   -x_1 x_2 \\
 &  &  x_1^2 +   x_2^2
  \end{array}\!\! \right]=x_0^4 x_1^2+x_1^4 x_0^2+x_2^4 x_0^2+x_1^4 x_2^2$$ 
  is a positive semidefinite quadratic matrix polynomial of a positive map that is not completely positive, but its determinant is an SOS polynomial.

Positive linear maps appear in matrix theory and operator algebras~\cite{P, Stormer}. For nearly half a century  (completely) positive maps gained a lot of attention in mathematical physics and lately in quantum computation and quantum information theory~\cite{Chuang, Cramer}. In the last decade completely positive maps are effectively used in semidefinite programming~\cite{KlepSchweighofer,  Lasserre}. In optimization completely positive maps are related to SOS matrix polynomials, which can be used to detect SOS-convex polynomials~\cite{Ahmadi}. In~\cite{HaMi1,HaMi2} the authors study nonnegative biquadratic forms arising from quasiconvex quadratic forms in relation with elasticity tensors.\\

\noindent
\textbf{Contributions}\hspace{3mm}
The main result of this paper is the construction of new families of positive linear maps $\Phi \colon \Sym_3\to \Sym_3$ such that the associated biquadratic forms $p_{\Phi}(\x,\y)=\y^T\Phi \left(\x\x^T \right) \y$ have 8, 9 or 10 zeros in $\mathbb{P}^2(\mathbb{R}) \times \mathbb{P}^2(\mathbb{R})$. Examples in Theorem~\ref{thm10zeros}, Theorem~\ref{thm9zeros} and Theorem~\ref{thm8zeros}  are the first known examples having 10, 9 and 8 zeros respectively. Note that the maximal posssible number of real zeros (if it is finite) is 10 and that such maps are never completely positive~\cite{Q}. 
From our construction, in which we prescribe the zero set of  biquadratic forms, we can directly deduce which nonnegative biquadratic forms are extremal. Our method of constructing positive maps by determining families of zeros of the associated biquadratic forms is to the best of our knowledge new. Moreover, we obtain significantly different examples as the ones found in the literature (see Remark~\ref{remzniclami}).

During our research we also learned that positive maps from $M_n(\mathbb{C})$ to $M_m(\mathbb{C})$ are of big importance in quantum computing. Positive maps that are not completely positive are entanglement witnesses of non-entanglement breaking maps~\cite{NieZhang}, whereby entanglement is most efficiently detected by the extremal positive maps. 
For this reason, in Section~\ref{SecExtension} we extend our examples of positive maps that are extreme rays in the cone of positive maps on $\Sym_3$  to positive maps that are extreme rays in the cone of positive maps on $M_ 3 (\mathbb{C})$.

Our explicit examples provide a big family of (extremal) positive maps on $\Sym_3$ that might help to understand the difference between the convex cones $ \Sigma_{B(3,3)} \subset P_{B(3,3)}$ and moreover their boundaries. Repeating the proofs of~\cite{Nie} and~\cite{Sturmf} for biquadratic forms, shows that the algebraic boundary of $P_{B(3,3)}$ is the \textit{discriminant} of biquadratic ternary forms. 
Using~\cite[Chapter 13. 2. B]{GelfKapr} we calculated that the discriminant is an irreducible hypersurface of degree 129 in $\mathbb{P}^{35}(\mathbb{R})$.\\

\noindent
\textbf{Acknowledgements}\hspace{3mm}
The authors would like to thank Bor Plestenjak for helping to numerically determine the regions of positivity long before we were able to provide proofs and to the referee for pointing out the difference between extremality of positive maps on $\Sym_3$ and  on $M_3(\mathbb{C})$. 

\section{Canonical form}
\label{SecCan}

In this section we construct a canonical form for positive maps under $\PGL_3 \times \PGL_3$  action.

If we write a positive linear map $\Phi \colon \Sym_3\to \Sym_3$ in coordinates
$$\Phi  \left( \left[
\begin{array}{ccc}z_{00}&z_{01}&z_{02}\\ &z_{11}&z_{12}\\ & &z_{22}
\end{array}
\right]\right) 
= z_{00} A_{00}+z_{01}A_{01}+z_{02}A_{02}+z_{11} A_{11}+z_{12}A_{12}+z_{22}A_{22},$$
where $A_{ij}$ are constant $3\times 3$ symmetric matrices
$$ A_{ij}=\left[
\begin{array}{ccc} a_{ij} & b_{ij} & c_{ij} \\
                              & d_{ij} & e_{ij}\\
                            &  & f_{ij}
\end{array}
\right],$$
then the corresponding nonnegative biquadratic form $p_{\Phi}$ on $\mathbb{P}^2\times \mathbb{P}^2$ equals
$$(y_0,y_1,y_2 ) . \left( 
 x_0^2 A_{00} + x_0 x_1 A_{01}+x_0 x_2 A_{02}+  
       x_1^2 A_{11}+ x_1 x_2 A_{12}+ x_2^2 A_{22} \right) .
        \left(\!\! \begin{array}{c} y_0 \\ y_1\\ y_2 \end{array}  \!\!\right). $$

Throughout the paper we will assume that $p_{\Phi}$ has at least 7 real zeros in $\mathbb{P}^2\times \mathbb{P}^2$. Moreover, we assume that 4 of the zeros have the following property: by projecting them to the first component we get 4 points in general position in $\mathbb{P}^2$, and the same holds by projecting them to the second component.
Since our aim is to construct positive maps that are not completely positive, these assumptions follow naturally from the following remark.
\begin{rmk} \label{remmaxzero}
 {\rm Quarez~\cite{Q} proved that a nonnegative biquadratic form which is not SOS can have at most 10 real zeros in $\mathbb{P}^2\times \mathbb{P}^2$. On the other hand, the number of real zeros of an SOS biquadratic form is either infinite or strictly less than 7. Nonnegative biquadratic forms with 7, 8, 9 or 10 real zeros therefore define positive maps that are not completely positive. \label{sosnuzeros}
In~\cite[Proposition 6.6]{Q} we also find that a biquadratic form has infinitely many zeros if any four zeros have the following property:  the projections of the four points to the first component are collinear in $ \mathbb{P}^2$, and the same holds for the projections to the second component. }
\end{rmk}

The group $\PGL_3 \times \PGL_3$ acts naturally on the convex cone of nonnegative biquadratic forms on $\mathbb{P}^2\times \mathbb{P}^2$ by 
$$ \y^T\Phi (\x \x^T)\y \mapsto (Q \y)^T  \Phi (P\x \x^T \!P^T)\, Q \y, $$
where $ P,Q \in \PGL_3$. Under this equivalence we can always assume that 
$$(1,1,-1;1,1,-1),\ (1,-1,1;1,-1,1),\ (-1,1,1;-1,1,1)\mbox{ and } (1,1,1;1,1,1)$$ are zeros of  $p_{\Phi}$. 
Since  $p_{\Phi}$ is nonnegative, it achieves its minimum on the real zero locus. Evaluating the gradient in 
the above 4 points gives a linear system of equations in the entries of $A_{ij}$,
$$\frac{\partial p_{\Phi}}{\partial x_i}=\frac{\partial p_{\Phi}}{\partial y_j}=0 \mbox{ for } i,j=0,1,2.$$
The solution is a $(36-20)$-parameter family
$$\begin{array}{l} \vspace{1mm}
A_{00}=\left[\begin{array}{ccc} a_{00} & b_{00} & c_{00} \\
                              & d_{00} & e_{00}\\
                               & & f_{00}
\end{array}  \right] \\ \vspace{1mm}
A_{01}= \left[\begin{array}{ccc}   -2 b_{00} &\frac{1}{2} \left(-a_{00}  - d_{11} + f_{22}- a_{11}- d_{00}  \right) & -e_{00} - e_{11} + e_{22} \\
        & -2 b_{11} & -c_{00} - c_{11} + c_{22}\\ 
           &  &  2 b_{22} 
       \end{array}  \right]\\ \vspace{1mm}
A_{11}= \left[\begin{array}{ccc} a_{11} & b_{11} & c_{11} \\
                            & d_{11} & e_{11}\\
                              &  & f_{11}
\end{array}  \right]\\ \vspace{1mm}
A_{02}=  \left[\begin{array}{ccc}  -2 c_{00} & -e_{00} + e_{11} - e_{22} & \frac{1}{2} (-a_{00}  + d_{11} - f_{22}- a_{22} -  f_{00} )\\ 
         & 2 c_{11} & -b_{00} + b_{11} - b_{22}\\
           & & -2 c_{22}
        \end{array}  \right]\\ \vspace{1mm}
 A_{12}= \left[\begin{array}{ccc}  2 e_{00} & c_{00} - c_{11} - c_{22} &  b_{00} - b_{11} - b_{22} \\
    & -2 e_{11}& \frac{1}{2} (a_{00}  - d_{11}  - f_{22}-d_{22}-  f_{11})\\ 
           & & -2 e_{22} 
             \end{array}  \right]    \\ \vspace{1mm}
 A_{22}= \left[\begin{array}{ccc} a_{22} & b_{22} & c_{22} \\
     & d_{22} & e_{22}\\ 
      &  & f_{22} \end{array}  \right]   ,
\end{array}$$
where $d_{00} + f_{00} =a_{11}+a_{22}$ and $a_{11} + f_{11}=d_{00}+d_{22}$. Note that the linearly dependent relation $a_{22}+d_{22}=f_{00}+f_{11}$ completes the symmetry. When a positive map is presented in the above basis, we say that it is in the \textit{canonical form}. 

Positive map $\Phi$ comes in pair with another positive map $\hat{\Phi}$, which is defined by the duality
\begin{equation} \y^T\Phi (\x \x^T)\y= \x^T \hat{\Phi} (\y \y^T)\x. \label{eq:dualhat} \end{equation} 
Note that $ \hat{\Phi}\colon \Sym_3 \to \Sym_3$ is exactly the adjoint map of $\Phi$. 
Clearly  $\Phi$ is completely positive if and only if  $\hat{\Phi}$ is completely positive.
For $\Phi$ in the canonical form, its pair $\hat{\Phi}\colon [w_{ij}]_{0 \leq i \leq j \leq 2} \to 
\sum_{0 \leq i \leq j \leq 2 } w_{ij} \hat{A}_{ij}$ equals 
$$ \begin{array}{l} \vspace{1mm}
\hat{A}_{00}=\left[\begin{array}{ccc} a_{00} & -b_{00} & -c_{00} \\
                              & a_{11} & e_{00}\\
                               & & a_{22}
\end{array}  \right] \\ \vspace{1mm}
\hat{A}_{01}= \left[\begin{array}{ccc}   2 b_{00} &\frac{1}{2} \left(-a_{00}  - d_{11} + f_{22}- a_{11}- d_{00}  \right) & -e_{00} +e_{11} - e_{22} \\
        & 2 b_{11} & c_{00} - c_{11} - c_{22}\\
           &  &  2 b_{22} 
       \end{array}  \right]\\ \vspace{1mm}
\hat{A}_{11}= \left[\begin{array}{ccc} d_{00} &-b_{11} & c_{11} \\
                            & d_{11} & -e_{11}\\
                              &  & d_{22}
\end{array}  \right]\\ \vspace{1mm}
\hat{A}_{02}=  \left[\begin{array}{ccc}  2 c_{00} & -e_{00} - e_{11} + e_{22} & \frac{1}{2} (-a_{00}  + d_{11} - f_{22}- a_{22} -  f_{00} )\\ 
         & 2 c_{11} & b_{00} - b_{11} - b_{22}\\
           & & 2 c_{22}
        \end{array}  \right]\\ \vspace{1mm}
 \hat{A}_{12}= \left[\begin{array}{ccc}  2 e_{00} & -c_{00} - c_{11} + c_{22} &  -b_{00} + b_{11} - b_{22} \\
    & 2 e_{11}& \frac{1}{2} (a_{00}  - d_{11}  - f_{22}-d_{22}-  f_{11})\\ 
           & & 2 e_{22} 
             \end{array}  \right]    \\ \vspace{1mm}
 \hat{A}_{22}= \left[\begin{array}{ccc} f_{00} & b_{22} & -c_{22} \\
     & f_{11} & -e_{22}\\ 
      &  & f_{22} \end{array}  \right],
\end{array}$$
where $d_{00} + f_{00} =a_{11}+a_{22},\ a_{11} + f_{11}=d_{00}+d_{22}$ and $a_{22}+d_{22}=f_{00}+f_{11}$ ensure that $\hat{\Phi}$ is also in the canonical form.

\begin{rmk}\label{remzniclami} {\rm The examples that we found in~\cite{CKL, CK, H,  LF, O, Q, TT} all satisfy 
$b_{ii}=c_{ii}=e_{ii}=0$ for $i=0,1,2$.}
\end{rmk}

\begin{rmk}\label{remkernel} {\rm 
Each additional  zero $T_k=(s_{0k},s_{1k},s_{2k};t_{0k},t_{1k},t_{2k})$ from the set of real zeros  of $ p_{\Phi}$
$$\{(\!1,\!1,\!-1;1,\!1,\!-1\!), (\!1,\!-1,\!1;1,\!-1,\!1\!), (\!-1,\!1,\!1;-1,\!1,\!1\!), (\!1,\!1,\!1;1,\!1,\!1\!)\} \cup \{ T_1,T_2,\ldots \},$$ 
imposes another set of linear equations in $a_{ii},b_{ii},c_{ii},d_{ii},e_{ii},f_{ii}$,
$$\frac{\partial p_{\Phi}}{\partial x_i}(T_k)=\frac{\partial p_{\Phi}}{\partial y_j}(T_k)=0 \mbox{ for } i,j=0,1,2.$$
We can see that, in order to obtain a nontrivial solution, increasing the number of real zeros also increases restrictions on the coordinates of zeros.  In practice this means that the two sets of points $\{(1,1,-1),(1,-1,1),(-1,1,1),(1,1,1),\} \cup \{(s_{0k},s_{1k},s_{2k})\}_k$ and $\{(1,1,-1),(1,-1,1),(-1,1,1),(1,1,1),\} \cup \{(t_{0k},t_{1k},t_{2k})\}_k$  
 are in some special geometric configuration in $\mathbb{P}^2$.
We believe that the reason for this hidden geometry lies in the following duality. For each zero $(s_0,s_1,s_2;t_0,t_1,t_2)$ of  $p_{\Phi}$, determinantal representation $\Phi (\x \x^T) $ is singular at  $(s_0,s_1,s_2)$ with eigenvector $(t_0,t_1,t_2)$, and  
 determinantal representation $\hat{\Phi} (\y \y^T)$ is singular at  $(t_0,t_1,t_2)$ with eigenvector $(s_0,s_1,s_2)$.
} \end{rmk}

\section{Nonnegative biquadratic forms with 10 zeros}
\label{sec:10zero}

For $t\in \mathbb{R}$ we will determine nonnegative biquadratic forms 
$p_t(\x,\y)=\y^T\Phi_t(\x\x^T)\y$ with the following real zeros,

$$\begin{array}{rcl} \vspace{1mm} 
\mathcal{Z}(p_t)&=& \left\{  (1,1,1;1,1,1), (1,1,-1;1,1,-1), (1,-1,1;1,-1,1), (-1,1,1;-1,1,1), \right. \\ \vspace{1mm}
 &&\ \, (1,t,0; t,1,0), (0,1,t;0,t,1), (t,0,1;1,0,t), \\ 
& &\  \, \left. (1,-t,0;-t,1,0), (0,1,-t;0,-t,1), (-t,0,1;1,0,-t)\right\}.
\end{array}$$

We demand that $p_t(\x,\y)$ is nonnegative, thus it is minimal on its real zero locus. This means that all partial derivatives evaluated on $\mathcal{Z}(p_t)$ must be 0.
The condition on the first four points puts $\Phi_t$ into the canonical form described in Section~\ref{SecCan}.
The remaining six points impose a system of equations that is linear in $a_{ii},b_{ii},c_{ii},d_{ii},e_{ii},f_{ii}$, and cubic in $t$ with the following nontrivial solutions:
\begin{eqnarray}
a_{00}=b_{00}=b_{11}=c_{00}=c_{22}=d_{11}=e_{11}=e_{22}=f_{22}=0,\label{eq:first11}  \\
a_{22}=f_{00}, \ d_{00}=a_{11},\  f_{11}=d_{22},\ t= \pm1; \nonumber \\
b_{00}=b_{11}=b_{22}=c_{00}=c_{11}=c_{22}=e_{00}=e_{11}=e_{22}=0, \label{eq:first1} \\
 a_{22}=d_{00}=f_{11}=0,\ f_{00}=a_{11}=d_{22},\ t=0; \nonumber \\
 b_{00}=b_{11}=b_{22}=c_{00}=c_{11}= c_{22}=e_{00}=e_{11} =e_{22}=0, \label{eq:genert} \\
 a_{00}=d_{11}=f_{22},\ 
f_{00}=a_{11}=d_{22}= \frac{a_{00}}{(-1 + t^2)^2},\  d_{00}=f_{11}=a_{22}= \frac{a_{00} t^4}{(-1 + t^2)^2}. \nonumber
\end{eqnarray}

It remains to be verified which of the above solutions define positive maps.
In the first case~\eqref{eq:first11} with $t=\pm 1$, the matrix
$$\Phi_{\pm 1}(\x \x^T)\!=\!\!
\left[\!\! \begin{array}{ccc}
\substack{a_{11} x_1^2 + 2 e_{00} x_1 x_2 + f_{00} x_2^2} & \substack{b_{22} x_2^2-a_{11} x_0 x_1 -e_{00} x_0 x_2 - c_{11} x_1 x_2}  &  \substack{c_{11} x_1^2 - e_{00} x_0 x_1-f_{00} x_0 x_2 - b_{22} x_1 x_2} \\
  &  \substack{a_{11} x_0^2 +2 c_{11} x_0 x_2+ d_{22} x_2^2} &  \substack{e_{00} x_0^2 - c_{11} x_0 x_1 - b_{22} x_0 x_2 - d_{22} x_1 x_2} \\ 
  & &  \substack{f_{00} x_0^2 + 2 b_{22} x_0 x_1 + d_{22} x_1^2}
   \end{array}\!\! \right]$$
is positive semidefinite (denoted by $\succeq 0$) iff all its principal minors are nonnegative for all $\x\in \mathbb{P}^2(\mathbb{R})$. Since the determinant is zero, the biquadratic form $p_{\pm 1}$ has infinitely many zeros and we obtain a completely positive map $\Phi_{\pm 1}$. The $1\times 1$ principal minors $\Phi_{\pm 1}(\x \x^T)$ are everywhere nonnegative when
 $$a_{11}\geq 0,\  d_{22}\geq 0 ,\  f_{00}\geq 0,\  a_{11} f_{00} -e_{00}^2 \geq 0 ,\  a_{11} d_{22} -c_{11}^2 \geq 0 ,\  d_{22} f_{00} -b_{22}^2 \geq 0,$$
moreover,  the three $2 \times 2$ principal minors of $\Phi_{\pm 1}(\x \x^T)$ are  everywhere nonnegative when
 $$\left[\! \begin{array}{ccc}
 -e_{00}^2 + a_{11} f_{00} & a_{11} b_{22} + c_{11} e_{00} & b_{22} e_{00} + c_{11} f_{00} \\
   & -c_{11}^2 + a_{11} d_{22} & b_{22} c_{11} + d_{22} e_{00} \\
 & & -b_{22}^2 + d_{22} f_{00}
  \end{array}\! \right] \succeq 0.$$
A straightforward calculation shows that these conditions are equivalent to 
$$\left[\! \begin{array}{ccc}
a_{11} & -c_{11} & -e_{00} \\
   & d_{22} & -b_{22} \\
 & &  f_{00}
  \end{array}\! \right] \succeq 0.$$
 
In the second case~\eqref{eq:first1} when $t=0$, the zero set $\mathcal{Z}(p_0)$ consists of 7 zeros. The positive semidefiniteness of
$$ \Phi_{0}(\x \x^T)\!=\!\!
\left[\!\! \begin{array}{ccc}
a_{00} x_0^2 +  a_{11} x_1^2 & -\frac{1}{2} ( a_{11} +a_{00} + d_{11} - f_{22}) x_0 x_1 & -\frac{1}{2} ( a_{11} +a_{00}  - d_{11} + f_{22}) x_0 x_2\\
& d_{11} x_1^2 + a_{11} x_2^2 &  -\frac{1}{2} ( a_{11} -a_{00}  + d_{11} + f_{22}) x_1 x_2 \\
 &  & a_{11} x_0^2 + f_{22} x_2^2 
  \end{array}\!\! \right]$$
can be again analyzed by the principal minors, which will be done in Section~\ref{sec:7zero}. 
  
In the remaining case \eqref{eq:genert} we take $a_{00}=(t^2-1)^2\neq 0$ and get
a new 1-parameter family of positive maps that are not completely positive.
\begin{thm} \label{thm10zeros} Linear maps $\Phi _t$ on $\Sym_3$ of the form 
$$\begin{array}{c}
\left[
\begin{array}{ccc} z_{00}& z_{01}&z_{02} \\
&   z_{11} &z_{12}\\
& &  z_{22}
\end{array}
\right]\\
\begin{tikzcd}
  \arrow[mapsto]{d}
  \\
\phantom{a} 
\end{tikzcd}\\
\left[
\begin{array}{ccc}(t^2-1)^2 z_{00}+z_{11}+t^4 z_{22}&-(t^4-t^2+1)z_{01}&-(t^4-t^2+1) z_{02} \\
&  t^4 z_{00}+(t^2-1)^2 z_{11}+z_{22} &-(t^4-t^2+1) z_{12}\\
& & z_{00}+t^4 z_{11}+ (t^2-1)^2 z_{22}
\end{array}
\right]\\
\end{array}$$  
are  positive for $ t\in \mathbb{R}.$
Apart from $t=\pm 1$, these positive maps are extremal and not completely positive. 
\end{thm}

\begin{proof}
The determinant of 
$$\Phi _t(\x \x^T)=\left[
\begin{array}{ccc}(t^2-1)^2x_0^2+x_1^2+t^4x_2^2&-(t^4-t^2+1)x_0 x_1&-(t^4-t^2+1)x_0 x_2 \\
&  t^4x_0^2+(t^2-1)^2x_1^2+x_2^2 &-(t^4-t^2+1)x_1 x_2\\
& & x_0^2+t^4x_1^2+ (t^2-1)^2x_2^2
\end{array}
\right]$$
is $(t^2-1)^2$ multiplied by the generalized Robinson polynomial
\begin{eqnarray}\label{eq:genrob}
t^4(x_0^6+x_1^6+x_2^6)-3(t^8-2t^6+t^4-2t^2+1)x_0^2x_1^2x_2^2+\\
+(t^8-2t^2)(x_0^4x_1^2+x_1^4x_2^2+x_2^4x_0^2)+(1-2t^6)(x_0^2x_1^4+x_1^2x_2^4+x_2^2x_0^4).\nonumber
\end{eqnarray}
This polynomial is nonnegative by~\cite[Example 6.5]{R}. 
The trace of the matrix $\Phi _t({{\tt x}}{{\tt x}}^T)$ equals 
$$\mathrm{Tr}(\Phi _t({\tt x}{\tt x^T}))=2(t^4-t^2+1)(x_0^2+x_1^2+x_2^2),$$
which is positive for each ${{\tt x}}\in \mathbb{P}^2(\mathbb{R})$. Similarly, the sum of all main $2\times 2$ subdeterminants of $\Phi _t({{\tt x}}{{\tt x}}^T)$ is 
$$(t^4-t^2+1)^2(x_0^2+x_1^2+x_2^2)^2,$$
which is again positive for all ${{\tt x}}$. 
This shows that  $\Phi _t$ is a positive linear map for all $t\in\mathbb{R}\backslash \{ \pm 1\}$.
   
Note that, if we substitute $t=0$ in the above matrix, we get (after a change of variables) the famous example of Choi considered in~\cite{Q} which is known to be extremal~\cite{H1}; and for $t=\pm 1$ we get an SOS matrix
$$\left[\!\!
\begin{array}{ccc} x_1^2+ x_2^2&-x_0 x_1&-x_0 x_2 \\
&  x_0^2+x_2^2 &-x_1 x_2\\
& &x_0^2+x_1^2
\end{array}\!\! \right]=
\left[\!\!\!
\begin{array}{c} x_1 \\
-\!x_0\\
0
\end{array}\!\!\! \right] [x_1 -\!\!x_0\, 0]+\left[\!\!\!
\begin{array}{c} -\!x_2 \\
0\\
x_0
\end{array}\!\!\! \right] [-\!x_2\, 0\, x_0]+\left[\!\!\!
\begin{array}{c} 0 \\
x_2\\
-\!x_1
\end{array}\!\!\! \right] [0\, x_2 -\!\!x_1].$$ 
An easy computation in \texttt{Mathematica} checks that for all other $t$ the generalized Robinson polynomial has 10 real zeros, therefore it is not a sum of squares by~\cite[Theorem 4.6]{CLR}, and hence $\Phi_t$ is not completely positive. 
It is obvious from our construction that $\Phi_t$ is extremal. Indeed, if  $\Phi_t=\Phi'_{t}+\Phi''_{t}$ would be a sum of two different positive maps, the associated biquadratic forms $p'_{t},  p''_{t}$ and $p_t$  would need to have the same 10 zeros. 
\end{proof}
\begin{cor}
As discussed in Remark~\ref{remmaxzero}, the upper bound on the number of real zeros of nonnegative biquadratic forms is 10. Theorem~\ref{thm10zeros} shows that the upper bound is obtained.
\end{cor}
\begin{rmk} {\rm Note that the generalized Robinson polynomial~\eqref{eq:genrob} converges to Robinson polynomial~\eqref{eq:therob} when $t \rightarrow \pm 1$. However,
$$\lim_{t \rightarrow \pm1} \det \Phi_{t}(\x \x^T)=0,$$
therefore the question whether there exists a positive semidefinite quadratic determinantal representation of Robinson polynomial is still open.
}\end{rmk}

\section{Nonnegative biquadratic forms with 9 zeros}
\label{sec:9zero}

For $p,q \in \mathbb{R}$ we will determine nonnegative biquadratic forms 
$p_{p,q}(\x,\y)=\y^T\Phi_{p,q}(\x\x^T)\y$ with the following zeros,
$$\begin{array}{l}\vspace{1mm}
 \left\{  (1,1,1;1,1,1), (1,1,-1;1,1,-1), (1,-1,1;1,-1,1), (-1,1,1;-1,1,1), \right. \\
\  \, \left.  (1,p,0; q,1,0), (1,-p,0;-q,1,0), (0,1,q;0,p,1), (0,1,-q;0,-p,1), (0,0,1;1,0,0)\right\}.
\end{array}$$
The same argument as in the previous section on the first four zeros, ensures that $\Phi_{p,q}$ is in the canonical form.
The condition that the partial derivatives of $p_{p,q}(\x,\y)$ are zero on the prescribed set of real zeros gives a linear system of equations in the coefficients of $p_{p,q}$ with the following nontrivial solutions
\begin{eqnarray}
b_{00}=b_{11}=b_{22}=c_{00}=c_{11}=c_{22}=e_{00}=e_{11}=e_{22}=0,\label{eq:first2}\\
a_{22}=d_{00}=f_{11}=0,\    a_{11}=d_{22}=f_{00},\  p=q=0; \nonumber \\
b_{00}=b_{11}=b_{22}=c_{00}=c_{11}=c_{22}=e_{00}=e_{11}=e_{22}=0, \label{eq:second2}\\
a_{00}=a_{22}=d_{00}=f_{11}=f_{22}=0,    \  a_{11}=d_{22}=f_{00}=-d_{11},\ p=0; \nonumber \\
 b_{00}=b_{11}=b_{22}=c_{00}=c_{11}=c_{22}=d_{00}=d_{11}=d_{22}=e_{00}=e_{11}=e_{22}=0, \label{eq:third2} \\
a_{11}=a_{22}= f_{00}=f_{11}=0, \  f_{22}=a_{00},\ q=0; \nonumber \\
 a_{22}=b_{22}=c_{00}=c_{22}=e_{00}=f_{00}=0, \  d_{00}=a_{11}, \  f_{11}=d_{22}, \label{eq:fourth2}\\
 a_{00}=f_{22}=\frac{d_{11}}{4},\  b_{00}=\frac{b_{11}}{2}, \
e_{22}=\frac{e_{11}}{2},\  p=\pm \frac{1}{\sqrt{2}}, \ q=\pm \sqrt{2}; \nonumber \\
a_{00}=a_{22}=b_{00}=b_{11}=b_{22}=c_{00}=c_{22}=e_{00}=e_{11}=e_{22}=f_{00}=f_{22}=0, \label{eq:fifth2}\\
d_{11}= 0,\ d_{00}=a_{11},\ d_{22}=f_{11},\ p=\frac{1}{q}; \nonumber \\
b_{00}=b_{11}=b_{22}=c_{00}=c_{11}=c_{22}=e_{00}=e_{11}=e_{22}=0,  \label{eq:sixth2}\\
a_{11}=a_{22}=d_{00}=d_{22}=f_{00}=f_{11}=0,\    a_{00}=f_{22}=\frac{d_{11}}{4},\ p=\frac{q}{2}; \nonumber \\
b_{00}=b_{11}=b_{22}=c_{00}=c_{11}=c_{22}=e_{00}=e_{11}=e_{22}=0, \ f_{00} = -\frac{ (2 p - q) (1 + p q)}{q ( p q-1)} d_{11}, \label{eq:seventh2}\\
a_{00}=f_{22}=\frac{p^2}{q^2}d_{11},\ a_{22}=0, \
d_{00}=f_{11}=\frac{p^2 (2 p - q) q}{( p q-1)^2} d_{11},  \
a_{11} =d_{22}=\frac{ 2 p - q}{q ( p q-1)^2} d_{11}. \nonumber 
\end{eqnarray}
Case~\eqref{eq:first2} with $p=q=0$ is the same as case~\eqref{eq:first1} with $t=0$, where we get the Choi set of 7 zeros which will be considered in Section~\ref{sec:7zero}.
In~\eqref{eq:second2} we get 
$$\Phi_{0,q}(\x\x^T)=d_{11}
\left[\!\! \begin{array}{ccc}
 -x_1^2  &0 & x_0 x_2\\
&  x_1^2 -x_2^2 &   0 \\
 &  &-x_0^2 
  \end{array}\!\! \right]$$ which is clearly not semidefinite.  On the other hand, case \eqref{eq:third2}  is SOS for $a_{00}>0$, since
 $$\Phi_{p,0}(\x\x^T)=a_{00}
\left[\!\! \begin{array}{ccc}
x_0^2  &0 & -x_0 x_2\\
& 0 &   0 \\
 &  &x_2^2 
  \end{array}\!\! \right]=a_{00}\left[\!\! \begin{array}{c} x_0 \\
0\\
-\!x_2
\end{array}\!\!\! \right] [x_0\, 0\, -\!x_2].$$

Determinantal representation $\Phi_{\substack{\substack{\pm \frac{1}{\sqrt{2} }, \pm \sqrt{2}}}}(\x\x^T)$ in~\eqref{eq:fourth2} equals 
$$\left[\!\!\!\! \begin{array}{ccc}
\substack{ \frac{d_{11}}{4}x_0^2 -b_{11} x_0 x_1+a_{11} x_1^2} \!\!\!&\!\!
\substack{ b_{11}\! (\!\frac{x_0^2}{2} \!+\! x_1^2\!) -(\!a_{11}\!+\! \frac{d_{11}}{2}\!)x_0 x_1 +\frac{e_{11}}{2} x_0 x_2-c_{11} x_1 x_2}  \!\!\!&\!\!
\substack{ c_{11} x_1^2-\frac{1}{2} e_{11} x_0 x_1+\frac{d_{11}}{4} x_0 x_2-\frac{b_{11}}{2} x_1 x_2} \\
 \!\!\!&\!\!\!  \substack{ a_{11} x_0^2 + d_{11} x_1^2 +  d_{22} x_2^2 - 2 (\!b_{11}\! x_0\! x_1 \! -\! c_{11}\! x_0\! x_2 \!+\! e_{11}\! x_1\! x_2\!) }  \!\!\!&\!  
 \substack{ e_{11}\! (\! x_1^2 \!+\!\frac{x_2^2}{2}\!)-c_{11} x_0 x_1+\frac{b_{11}}{2} x_0 x_2- (\!\frac{d_{11}}{2}\!+\!d_{22}\!) x_1 x_2 }\\
  \!\!\!&\!\!\!    \!\!\!&\! \substack{ d_{22} x_1^2 -e_{11} x_1 x_2+ \frac{d_{11}}{4} x_2^2 }
  \end{array}\!\!\!\! \right]$$
  with determinant
$$\frac{1}{4} (-c_{11}^2 d_{11} - b_{11}^2 d_{22} - a_{11} e_{11}^2+ a_{11} d_{11} d_{22} + 2 b_{11} c_{11} e_{11} ) x_1^2 (x_0^2 - 2 x_1^2 + x_2^2)^2.$$
This means that  $p_{\substack{\substack{\pm \frac{1}{\sqrt{2} }, \pm \sqrt{2}}}}(\x,\y)$ has infinitely many real zeros and is SOS whenever 
$$\Phi_{\substack{\substack{\pm \frac{1}{\sqrt{2} }, \pm \sqrt{2}}}}(\x\x^T)\succeq 0\ \mbox{ for all }\ \x\in 
\mathbb{P}^2(\mathbb{R}).$$  
We invite the reader to prove that the principal minors of  $\Phi_{\substack{\substack{\pm \frac{1}{\sqrt{2} }, \pm \sqrt{2}}}}(\x\x^T)$ 
are everywhere nonnegative if and only if
$$ \left[ \begin{array}{ccc}
a_{11} & -b_{11} & c_{11} \\
 & d_{11} & -e_{11} \\
 & & d_{22} 
  \end{array} \right] \succeq 0.$$
(This calculation is a good exercise in semidefiniteness of binary and ternary quadratic forms that will feature in many other examples). 
  
Both cases \eqref{eq:fifth2} and \eqref{eq:sixth2} 
$$\Phi_{p,\frac{1}{p}}(\x\x^T)=
\left[\!\!\! \begin{array}{ccc}
a_{11} x_1^2 & -x_1 (a_{11} x_0 + c_{11} x_2) & c_{11} x_1^2 \\
 & a_{11} x_0^2 +  2 c_{11} x_0 x_2 + d_{22} x_2^2 & -x_1 (c_{11} x_0 + d_{22} x_2)\\
 & & d_{22} x_1^2
 \end{array}\!\!\! \right]$$
 and
   $$\Phi_{p,2p}(\x\x^T)=\frac{d_{11}}{4}
   \left[\!\!\! \begin{array}{ccc}
   x_0^2 & -2 x_0 x_1 &  x_0 x_2 \\
   &  4 x_1^2 &  -2 x_1 x_2 \\
   & &  x_2^2
 \end{array}\!\!\! \right] = \frac{d_{11}}{4}
 \left[\!\! \begin{array}{c} x_0 \\
-\! 2 x_1\\
x_2
\end{array}\!\!\! \right] [x_0\, -\! 2 x_1\, x_2]$$  
have determinant 0. This implies that  $\Phi_{p,\frac{1}{p}}$ is completely positive (since the associated biquadratic form has infinitely many real zeros) when $a_{11} \geq 0,\, d_{22}\geq 0$ and $a_{11}  d_{22}-c_{11}^2\geq 0$. 
Indeed, under these conditions  $\Phi_{p,\frac{1}{p}}(\x\x^T)$ is positive semidefinite for all 
$\x \in \mathbb{P}^2(\mathbb{R})$. Obviously, $\Phi_{p,2p}$ is  completely positive for $d_{11}\geq 0$.
  
The most interesting case is~\eqref{eq:seventh2}, where  we get a new family 
of positive maps that are not completely positive.
\begin{thm} \label{thm9zeros} Linear maps on $\Sym_3$ of the form 
$$
\begin{array}{c}
\left[
\begin{array}{ccc} z_{00}& z_{01}&z_{02} \\
&   z_{11} &z_{12}\\
& &  z_{22}
\end{array}
\right]\\
\begin{tikzcd}
  \arrow[mapsto]{d}
  \\
\phantom{a} 
\end{tikzcd}\\
\! \left[\!\!\! \begin{array}{ccc}
 \substack{ p^2(p q-1)^2} z_{00} + \substack{q(2 p - q) }z_{11} \!\!\!\!\!\!\! & \!\!\!\!\!\!\! \substack{ -p q (1 -  q^2 + p^2 q^2) }z_{01}  
 \!\!\!\!\!\!\! & \!\!\!\!\!\!\! 
\substack{ (p q-1)(p^2 + p q - p^3 q - q^2 + p^2 q^2)} z_{02}   \\
\!\!\!\!\!\!\! & \!\!\!\!\!\!\!  \substack{ p^2 q^3 (2 p - q)}  z_{00} \!+\! \substack{q^2 ( p q-1)^2} z_{11} \!+\! \substack{ q(2 p-q) } z_{22}    \!\!\!\!\!\!\! & \!\!\!\!\!\!\!  \substack{ -p q (1 -  q^2 + p^2 q^2)} z_{12}  \\
   \!\!\!\!\!\!\! & \!\!\!\!\!\!\!  \!\!\!\!\!\!\!\! & \!\!\!\!\!\!\!\!   \substack{q (2 p - q) (1 - p^2 q^2) } z_{00} \!+\! \substack{p^2 q^3 (2 p - q) } z_{11} \!+\! \substack{ p^2(p q-1)^2 }z_{22}
 \end{array}\!\!\! \right]
 \end{array}$$  
are  positive for 
$$\left\{ (p,q)\in [0,1/\sqrt{2}] \times [0,\sqrt{2}] \colon 2 p - q \geq 0,\ (p^2 - 1)^2 q^2 - p^2\geq 0 \right\}.$$
For all $(p,q)$ in this region, except on the line $q=2p$, these positive maps are extremal and not completely positive. 
\end{thm}
\begin{proof}
 If we take $d_{11}=q^2 (p q-1)^2$ in~\eqref{eq:seventh2}, we get that $\Phi_{p,q}(\x\x^T)$ equals
 $$\! \left[\!\!\! \begin{array}{ccc}
 \substack{ p^2(p q-1)^2 x_0^2 + q(2 p - q) x_1^2}   & \substack{ -p q (1 -  q^2 + p^2 q^2) x_0 x_1} &
\substack{ (p q-1)(p^2 + p q - p^3 q - q^2 + p^2 q^2)x_0 x_2  } \\
   &\substack{ p^2 q^3 (2 p - q)  x_0^2 +q^2 ( p q-1)^2 x_1^2 + q(2 p-q) x_2^2 }&\substack{ -p q (1 -  q^2 + p^2 q^2) x_1 x_2}  \\
   & & \substack{q (2 p - q) (1 - p^2 q^2)  x_0^2 + p^2 q^3 (2 p - q)  x_1^2+ p^2(p q-1)^2 x_2^2}
 \end{array}\!\!\! \right]\!\!.$$  
 Observe that its first and second leading principal minors 
 \begin{eqnarray*}
 p^2(p q-1)^2 x_0^2 + q(2 p - q) x_1^2\  \mbox{ and}\\
 (2 p - q) q  \left( q^2 (-1 + p q)^2 (p x_0 - x_1)^2 (p x_0 + x_1)^2 + (p^2 (1 - p q)^2 x_0^2 +
     q (2 p - q) x_1^2) x_2^2 \right) \hspace{1cm}
 \end{eqnarray*}
 are nonnegative for all $\x$ if and only if $q(2 p - q)\geq 0$,  and that  its 
determinant equals
$q^2( p q-1)^2 (2 p - q)^2$ times 
\begin{eqnarray*} p^4 q^2 (1 - p^2 q^2) x_0^6 + p^2 q^2 (-2 + 2 p^2 q^2 + p^4 q^2) x_0^4 x_1^2+ p^2 (1 - 2 p^2 q^2 + q^4 - 2 p^2 q^4 + p^4 q^4) x_0^4 x_2^2\\ 
+q^2 (1 -  p^2 q^2 - 2 p^4 q^2) x_0^2 x_1^4- ( p^2 + 2 q^2 - 6 p^2 q^2 + p^2 q^4 - 
  2 p^4 q^4 + p^6 q^4) x_0^2 x_1^2 x_2^2 \\
  +(q^2 -p^2+ p^4 q^2 -  2 p^2 q^2) 
x_0^2 x_2^4
+ p^2 q^4  x_1^6-2 p^2 q^2 x_1^4 x_2^2+p^2 x_1^2 x_2^4.
\end{eqnarray*}
First we consider the cases when the above determinant is constantly 0:
\begin{itemize}
\item
for $q=0$ we get the congruence map $\Phi_{p,0}$ obtained in \eqref{eq:third2} with $a_{00}=p^2$; 
 \item
for $q=\frac{1}{p}$ we get the completely positive map $\Phi_{p,\frac{1}{p}}$ obtained in \eqref{eq:fifth2}  with
$c_{11}=0$ and $a_{11}=d_{22}=\frac{2 p^2-1}{p^2}\geq 0$, which in particular implies $|p| \geq \frac{1}{\sqrt{2}}$; 
\item 
for  $q=2p$ we get the congruence map $\Phi_{p,2p}$ obtained in \eqref{eq:sixth2} with  $d_{11}=4 p^2 (1 - 2 p^2)^2$.
 \end{itemize}
From now on we assume that $q ( p q-1)  (2 p - q) \neq 0$. We will determine the conditions on $p,q$ under which  
$\Phi_{p,q}$ will be a positive linear map. Since the determinant is not constantly 0 and the principal minors of  $\Phi_{p,q}(\x\x^T)$ are continuous in $\x$, it is by Sylvester's criterion enough to verify that the leading principal minors are nonnegative for all 
$\x \in \mathbb{P}^2(\mathbb{R})$.

Note that the determinant is a quadric in $x_2^2$, thus it is everywhere nonnegative if and only if for all real $x_0,x_1$ one of the following holds:
 \begin{eqnarray*} 
\mbox{(i)} & \mbox{the coefficient at } x_2^4  \mbox{ is } & ((p^2 - 1)^2 q^2 - p^2) x_0^2 + p^2  x_1^2\geq 0, \\
&  \mbox{the constant term}   \mbox{ is } & q^2 (p^2 x_0^2 - x_1^2)^2 \left((1 - p^2 q^2) x_0^2 + p^2 q^2 x_1^2\right)\geq 0,
  \\
&  \mbox{and the discriminant  is }  & p^2 (1 - q^2 + p^2 q^2)^2 x_0^2 (x_0 - x_1)^2 (x_0 + x_1)^2 \times \\
& &(p^2 (-1 - q^2 + p^2 q^2)^2 x_0^2+  4 (p^2-1) q^2 x_1^2) \leq 0,
 \end{eqnarray*}
 or
 \begin{eqnarray*}
\mbox{(ii)} & \mbox{the coefficient at } x_2^4  \mbox{ is } & ((p^2 - 1)^2 q^2 - p^2) x_0^2 + p^2  x_1^2\geq 0, \\
& \mbox{the coefficient at } x_2^2 \mbox{ is } & p^2 ((1-p^2)^2 q^4- 2 p^2 q^2  +1) x_0^4 -\\
& & \!\!\!\!\!\!\!\!\!\!\!\!\!\!\!\!\!\!(p^2 + 2 q^2 - 6 p^2 q^2 + p^2 (1 - p^2)^2 q^4) x_0^2 x_1^2 -2 p^2 q^2 x_1^4 \geq 0,  \\
& \mbox{and the constant term is }  & q^2 (p^2 x_0^2 - x_1^2)^2 \left((1 - p^2 q^2) x_0^2 + p^2 q^2 x_1^2\right)\geq 0.
 \end{eqnarray*}
 We have already shown that the nonnegativity of the first two leading principal minors of $\Phi_{p,q}(\x\x^T)$  implies $q(2 p - q)\geq 0$.
 Moreover, nonnegativity of the leading coefficient and of the constant term in the determinant  implies $(p^2 - 1)^2 q^2 - p^2\geq 0$ and $1 - p^2 q^2\geq 0$.
 The  intersection of these regions 
 $$\mathcal{R}=\left\{ (p,q)\in \mathbb{R}^2 \colon q(2 p - q)\geq 0,\ (p^2 - 1)^2 q^2 - p^2\geq 0,\  1 - p^2 q^2\geq 0 \right\}$$ 
 is shown in Figure~\ref{pictRegion9pts}.
 \begin{figure}
\begin{center}
\includegraphics[width=8cm]{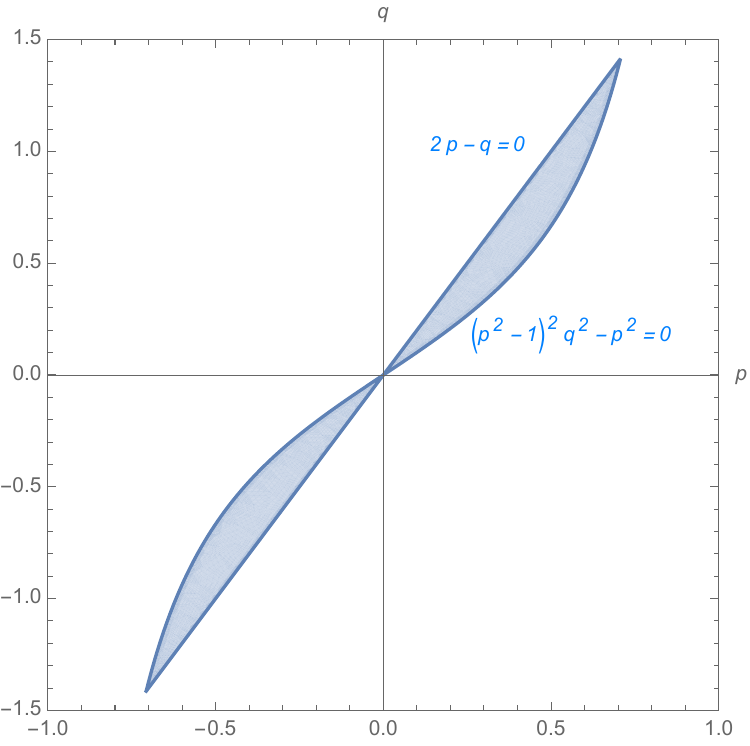}
\end{center}
\caption{Intersection of regions $q(2 p - q)\geq 0,\ (p^2 - 1)^2 q^2 - p^2\geq 0,\  1 - p^2 q^2\geq 0$. }
\label{pictRegion9pts} 
\end{figure}
A simple calculation shows that $(p,q)\in [-1/\sqrt{2},1/\sqrt{2}] \times [-\sqrt{2},\sqrt{2}]$ and that $\mathcal{R}$ is bounded by the curves with equations $2 p - q=0$ and $(p^2 - 1)^2 q^2 - p^2=0$. 
Since  $\Phi_{p,q}$  and the set of zeros of $p_{p,q}$ are both invariant under the simultaneous sign change of $\pm p,\pm q$, it is enough to consider $0\leq p\leq \frac{1}{\sqrt{2}} $ and $0\leq q \leq \sqrt{2}$. 

We can now prove that $\Phi_{p,q}(\x \x^T)\succeq 0$ holds for $(p,q)\in \mathcal{R}$ and all $\x \in \mathbb{P}^2(\mathbb{R})$.
For $$ p^2 (-1 - q^2 + p^2 q^2)^2 x_0^2  \leq 4 q^2 (1-p^2) x_1^2$$ the discriminant in (i) is clearly not positive. It remains to show that for $(p,q)\in \mathcal{R}$ and 
$ p^2 (-1 - q^2 + p^2 q^2)^2 x_0^2  \geq 4 q^2 (1-p^2) x_1^2$ the coefficient at $x_2^2$ in (ii) is nonnegative. 
View 
$$p^2 ( (1-p^2)^2 q^4- 2 p^2 q^2  +1) x_0^4 -
(p^2 + 2 q^2 - 6 p^2 q^2 + p^2 (1-p^2)^2 q^4) x_0^2 x_1^2 -2 p^2 q^2 x_1^4$$
as a quadric in $x_1^2$. Its leading term $-2 p^2 q^2$ is clearly negative for $(p,q)\neq(0,0)$ and the constant term 
$P=p^2  \left( \left(1-p^2\right)^2 q^4- 2 p^2 q^2  +1\right)$ is positive on $\mathcal{R} \backslash \{(1/\sqrt{2}, \sqrt{2})\}$ 
since
 $$\left(1-p^2\right)^2 q^4- 2 p^2 q^2  +1= \left(\left(p^2-1 \right)^2 q^2-p^2\right) q^2  +1-  p^2 q^2\geq 0,$$
as illustrated in Figure~\ref{RegionII9pts}. 
\begin{figure}
\begin{center}
\includegraphics[width=8cm]{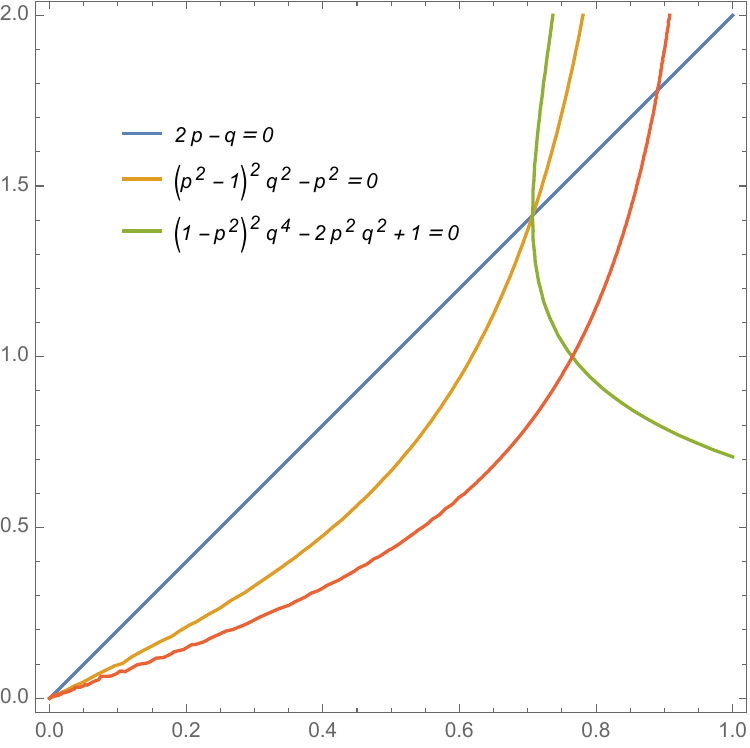}
\end{center}
\caption{Illustration of  nonnegativity of the coefficient at $x_2^2$ in (ii) on $\mathcal{R}$}
\label{RegionII9pts} 
\end{figure}
Denote 
$B= -p^2 - 2 q^2 + 6 p^2 q^2 - p^2 (1-p^2)^2 q^4$. Then the coefficient at $x_2^2$ in (ii) equals
\begin{eqnarray*}
P x_0^4 + B x_0^2 x_1^2 -2 p^2 q^2 x_1^4=\\
\left( x_1^2 + \frac{-B + \sqrt{B^2 + 8 P p^2 q^2} }{4 p^2 q^2}x_0^2  \right)
\left( -2 p^2 q^2 x_1^2 + \frac{B + \sqrt{B^2 + 8 P p^2 q^2} }{2}x_0^2  \right),
\end{eqnarray*}
which is nonnegative if and only if the second factor is nonnegative. From the assumption 
$ p^2 (-1 - q^2 + p^2 q^2)^2 x_0^2  \geq 4 q^2 (1-p^2) x_1^2$
it follows that
\begin{eqnarray*}
-2 p^2 q^2 x_1^2 + \frac{B + \sqrt{B^2 + 8 P p^2 q^2} }{2}x_0^2 \geq 
- p^4 \frac{ (-1 - q^2 + p^2 q^2)^2}{2  (1 - p^2)} x_0^2 + \frac{B + \sqrt{B^2 + 8 P p^2 q^2} }{2} x_0^2.
\end{eqnarray*}
We claim that the right hand side is nonnegative on $\mathcal{R}$, which is equivalent to
$$\sqrt{B^2 + 8 P p^2 q^2} \geq  p^4 \frac{ (-1 - q^2 + p^2 q^2)^2}{  (1 - p^2)}-B.$$
For the last inequality it is enough to show that it holds for the absolute value of the right hand side.
After squaring, we get 
 $$ \frac{p^4 (1 - q^2 + p^2 q^2)^2}{(1-p^2)^2}
 \left(2 q^2-p^2-3 p^2 q^2+p^4 q^2 \right) \left(2-p^2-p^2 q^2+p^4 q^2 \right) \geq 0,$$
 which holds  since the last two factors are both nonnegative on $\mathcal{R}$.
 Indeed, 
$$2 q^2-p^2-3 p^2 q^2+p^4 q^2=\left(p^2-1 \right)^2 q^2-p^2 +q^2(1-p^2)\geq 0$$
and
$$2-p^2-p^2 q^2+p^4 q^2=2-p^2+p^2 q^2(p^2-1)\geq 2-p^2+p^2-1=1.$$
The above manipulation with inequalities is illustrated in Figure~\ref{RegionII9pts}, where the red curve has equation
$$- p^4 \frac{ (-1 - q^2 + p^2 q^2)^2}{2  (1 - p^2)}  + \frac{B + \sqrt{B^2 + 8 P p^2 q^2} }{2}=0.$$ 
This concludes the proof of positivity of $\Phi_{p,q}$.

Finally we check for which $(p,q)\in \mathcal{R}$ the positive map $\Phi_{p,q}$ is completely positive or equivalently, when $p_{p,q}$ is SOS.  By Remark~\ref{remmaxzero} it suffices to prove that $p_{p,q}$ does not have infinitely many real zeros.  On the boundary $q^2=\frac{p^2}{(p^2-1)^2}$ of $\mathcal{R}$, the real zeros of  
$\det \Phi_{p,q}(\x \x^T)$ are the real zero set of 
\begin{eqnarray*}(p^4 - 2 p^6) x_0^6+ p^2 (-2 + 4 p^2 + p^6)x_0^4 x_1^2 -( p^2-1)^2 ( 2 p^2-1)x_0^4 x_2^2+(1 - 
 2 p^2 - 2 p^6)x_0^2 x_1^4  \\ -(p^2-1)^2 (3 - 8 p^2 + 
   2 p^4)x_0^2 x_1^2 x_2^2+ p^4 x_1^6 -2 p^2 (p^2-1)^2 x_1^4 x_2^2+ (p^2-1)^4  x_1^2 x_2^4,
   \end{eqnarray*}
which for $p\in (0,1/\sqrt{2})$ equals
$$\begin{array}{l}\vspace{1mm}
 \left\{  (1,1,1), (1,1,-1), (1,-1,1), (-1,1,1), \right. \\
\  \, \left.  (1,p,0), (1,-p,0), (0,1,\frac{p}{1-p^2}), (0,1,-\frac{p}{1-p^2}), (0,0,1)\right\}.
\end{array}$$  
The kernels of $\Phi_{p,q}(\x\x^T)$ in these 9 points are one dimensional, which means that $p_{p,q}$ has 9 real zeros and is thus not SOS.
On the other boundary $q=2p$ of $\mathcal{R}$ we get a congruence map already considered in case~\eqref{eq:sixth2}. 
For $(p,q)$  in the interior of  $\mathcal{R}$ we used  \texttt{Mathematica} to solve the system 
$$\frac{\partial \det \Phi_{p,q}(\x \x^T)}{\partial x_i}=0,\ i=0,1,2 $$
and obtained 9 solutions 
$$\left\{  (1,1,1), (1,1,-1), (1,-1,1), (-1,1,1), (1,p,0), (1,-p,0), (0,1,q), (0,1,-q), (0,0,1)\right\},$$
which are the real zeros of  $\det \Phi_{p,q}(\x \x^T)$.
Additionally we need to show that $\Phi_{p,q}(\x \x^T)$ has one dimensional kernels at these points. 
By Remark~\ref{remkernel} it is equivalent to show that $\det \hat{\Phi}_{p,q} (\y \y^T)$ has 9 real zeros, where $\hat{\Phi}_{p,q}$ is the dual map of $\Phi_{p,q}$ defined in~\eqref{eq:dualhat}. 
For the latter we again used  \texttt{Mathematica}  to prove that
$$\left\{  (1,1,1), (1,1,-1), (1,-1,1), (-1,1,1), (q,1,0), (-q,1,0), (0,p,1), (0,-p,1), (1,0,0)\right\}$$
are the only solutions  of the system
$$\frac{\partial \det \hat{\Phi}_{p,q}(\y \y^T)}{\partial y_i}=0,\ i=0,1,2 .$$
We proved that $p_{p,q}$ has 9 real zeros for $(p,q)$ in the interior of  $\mathcal{R}$, thus by Remark~\ref{remmaxzero} it is not SOS and $\Phi_{p,q}$ is not a completely positive map. 

From our construction it follows that for each  
$(p,q) \in \mathcal{R} \backslash \left\{ (0,0), \left(\frac{1}{\sqrt{2}},  \sqrt{2} \right) , \left( - \frac{1}{\sqrt{2}},  - \sqrt{2} \right)  \right\}$ there is up to a scalar a unique nonnegative biquadratic form $p_{p,q}$ with the prescribed 9 real zeros, therefore all such positive maps $\Phi_{p,q}$ are extremal. 
\end{proof}
\medskip

\section{Nonnegative biquadratic forms with 8 zeros}
\label{sec:8zero}

For $m,n \in \mathbb{R}$ we will determine nonnegative biquadratic forms 
$p_{m,n}(\x,\y)=\y^T\Phi_{m,n}(\x\x^T)\y$ with the following zeros,
$$\begin{array}{l}\vspace{1mm}
 \left\{  (1,1,1;1,1,1), (1,1,-1;1,1,-1), (1,-1,1;1,-1,1), (-1,1,1;-1,1,1), \right. \\
\  \, \left.  (1,0,0; m,1,0), (1,n,0; 0,1,0), (0,1,0;0,0,1), (0,0,1;1,0,0)\right\}.
\end{array}$$
The condition that the partial derivatives of $p_{m,n}(\x,\y)$ are zero on the set of zeros gives a linear system of equations in the coefficients of $p_{m,n}$ in the canonical form, with the following nontrivial solutions

\begin{eqnarray}
b_{00}=b_{11}=b_{22}=c_{00}=c_{11}=c_{22}=e_{00}=e_{11}=e_{22}=0, \label{eq:mnis0} \\
a_{22}= d_{00}=f_{11}=0, \  a_{11}=d_{22}=f_{00}, \ m=n=0; \nonumber \\
b_{00}=b_{11}=b_{22}=c_{00}=c_{11}=c_{22}=e_{00}=e_{11}=e_{22}=0, \label{eq:mis0}\\
a_{22}=d_{00}=d_{11}=f_{11}=0,    \  a_{11}=d_{22}=f_{00}, \ f_{22}=a_{00}+a_{11},\ m=0; \nonumber \\
b_{00}=b_{11}=b_{22}=c_{00}=c_{11}=c_{22}=e_{00}=e_{11}=e_{22}=0,\label{eq:nis0}\\
a_{00}=a_{22}=d_{00}=f_{11}=0,\    a_{11}=d_{22}=f_{00},\  f_{22}=a_{11}+d_{11},\ n=0; \nonumber \\
a_{22}=b_{22}=c_{00}=c_{11}=c_{22}=e_{00}=e_{11}=e_{22}=f_{11}=0,    \label{eq:mnnot0} \\
  a_{00} = \frac{d_{00}}{m^2},\  b_{00} =-\frac{d_{00}}{m}, \  b_{11} = \frac{d_{00}}{n}, \ 
  d_{11} = \frac{d_{00}}{n^2}, \nonumber\\
  d_{22}=f_{00}=a_{11}-d_{00} ,\ f_{22}=a_{11}-d_{00} \left(1 - \frac{1}{m^2} - \frac{1}{n^2} - \frac{2}{m n} \right) .\nonumber 
\end{eqnarray}

In \eqref{eq:mnis0} with $m=n=0$  the biquadratic form $p_{0,0}$ has the 7 zeros of Choi, which we will analyze in Section~\ref{sec:7zero}.
In~\eqref{eq:mis0} we get
$$\Phi_{0,n}(\x\x^T)=
\left[\!\! \begin{array}{ccc}
 a_{00} x_0^2 + a_{11} x_1^2 &0 & -(a_{00} + a_{11}) x_0 x_2\\
&    a_{11} x_2^2 & -a_{11} x_1 x_2   \\
 &  & a_{11} x_0^2 + (a_{00} + a_{11}) x_2^2
  \end{array}\!\! \right]$$ and $\det \Phi_{0,n}(\x\x^T)=
  a_{11}^2 (x_0^2 - x_1^2)  x_2^2 \left(a_{00} x_0^2 +a_{11} x_1^2 - \left(a_{00}+a_{11}\right) x_2^2 \right)$. Such $\Phi_{0,n}$
  is not a positive map unless $a_{11}=0$, in which case it is a completely positive congruence map.
Analogously, case \eqref{eq:nis0}  is not semidefinite  for $d_{00}\geq 0$ and  $a_{11}>0$, since
 $$\Phi_{m,0}(\x\x^T)=
\left[\!\! \begin{array}{ccc}
 a_{11} x_1^2 &0 & -a_{11} x_0 x_2\\
&  d_{11} x_1^2 + a_{11} x_2^2 & -(a_{11} + d_{11}) x_1 x_2   \\
 &  & a_{11} x_0^2 + (a_{11} + d_{11}) x_2^2
  \end{array}\!\! \right]$$ 
  has  determinant
$$  a_{11}^2 (x_1^2 - x_2^2) \left(d_{11} x_0^2 x_1^2 + a_{11} x_0^2 x_2^2 - 
   (a_{11}+d_{11}) x_1^2 x_2^2 \right),$$ and for $d_{00}\geq 0$ and  $a_{11}=0$ the map $\Phi_{m,0}$ is a congruence map.
   
The last case~\eqref{eq:mnnot0} yields a three dimensional family of positive maps which we describe in the following theorem.
 \begin{thm} \label{thm8zeros} Linear maps on $\Sym_3$ of the form 
$$\begin{array}{r} \vspace{2mm} 
\left[
\begin{array}{ccc} z_{00}& z_{01}&z_{02} \\
&   z_{11} &z_{12}\\
& &  z_{22}
\end{array}
\right] \hspace{5mm} \mapsto \hspace{2.5cm} 
b \left[\!\!\! \begin{array}{ccc}
 z_{11} & 0 & -z_{02} \\
 & z_{22} & - z_{12}  \\
 & & z_{00} + z_{22}
 \end{array}\!\!\! \right]+  \\ 
   \left[\!\!\! \begin{array}{ccc}
n^2 (z_{00} +2  m z_{01}+m^2 z_{11}) & - m n \left(n z_{00}+(m n-1)z_{01}-m z_{11}\right) & - n (m +
      n) (z_{02} + m z_{12})      \\
 &   m^2 (n^2 z_{00} -2 n z_{01}+z_{11}) &  m (m + n) (n z_{02} - z_{12})    \\
 & &  (m + n)^2 z_{22}
 \end{array}\!\!\! \right]
 \end{array}  $$
are  positive for 
$$\left\{ (m,n,b)\in [-1,1]^2 \times [0,1/4] \colon 0\leq b\leq \min \{-(2 m n \!+\! n^2 \!+\! m^2 n^2),-(2 m n \!+\!  m^2 n^2 \!+\! m^2)\} \right\}.$$
These positive maps are completely positive for $b=0$ and not completely positive elsewhere. Moreover, for  
$b= \min \left\{-(2 m n + n^2 + m^2 n^2),-(2 m n +  m^2 n^2 +m^2) \right\}$ they are extremal positive maps.
\end{thm}
\begin{proof}
If we multiply $\Phi_{m,n}(\x\x^T)$ in~\eqref{eq:mnnot0} by $m^2 n^2$, we get
\begin{eqnarray*}
a_{11} m^2 n^2 \left[\!\!\! \begin{array}{ccc}
 x_1^2 & 0 & -x_0 x_2\\
 & x_2^2 & - x_1 x_2 \\
 & & x_0^2 + x_2^2
 \end{array}\!\!\! \right]+\\
 d_{00}  \left[\!\!\! \begin{array}{ccc}
 \substack{n^2 (x_0^2 + 2 m x_0 x_1)} & \substack{-m n (n x_0 - x_1) (x_0 + m x_1)} & \substack{-n (m x_0 + 
     n x_0 - m^2 n x_0 + m^2 x_1 + m n x_1) x_2} \\
 &  \substack{m^2 (n^2 x_0^2 - 2 n x_0 x_1 + x_1^2 - n^2 x_2^2)} & 
   \substack{m (m n x_0 + n^2 x_0 - m x_1 - n x_1 + m n^2 x_1) x_2}  \\
 & &   \substack{-m^2 n^2 x_0^2 + 
   (m+n)^2 x_2^2 - m^2 n^2 x_2^2}
 \end{array}\!\!\! \right]=\\
 (a_{11}-d_{00}) m^2 n^2 \left[\!\!\! \begin{array}{ccc}
 x_1^2 & 0 & -x_0 x_2\\
 & x_2^2 & - x_1 x_2 \\
 & & x_0^2 + x_2^2
 \end{array}\!\!\! \right]+\\
 d_{00}  \left[\!\!\! \begin{array}{ccc}
n^2 (x_0 + m x_1)^2 & - m n (n x_0 - x_1) (x_0 + m x_1) & - n (m +
      n) (x_0 + m x_1) x_2     \\
 &   m^2 (n x_0 - x_1)^2 &  m (m + n) (n x_0 - x_1) x_2   \\
 & &  (m + n)^2 x_2^2
 \end{array}\!\!\! \right].
 \end{eqnarray*}  
Therefore we need to determine $m, n$ and $b$ for which $\Phi_{b,m,n}(\x \x^T)=$
\begin{eqnarray*}
b \left[\!\!\! \begin{array}{ccc}
 x_1^2 & 0 & -x_0 x_2\\
 & x_2^2 & - x_1 x_2 \\
 & & x_0^2 \!+\! x_2^2
 \end{array}\!\!\! \right]+
\left[\!\! \begin{array}{c} -n (x_0 \!+\! m x_1) \\
m (n x_0 \!-\! x_1) \\
(m \!+\! n) x_2
\end{array}\!\!\! \right] [-n (x_0 \!+\! m x_1), m (n x_0 \!-\! x_1), (m \!+\! n) x_2]=\\
 \left[\!\!\! \begin{array}{ccc}
b x_1^2+ n^2 (x_0 + m x_1)^2 & - m n (n x_0 - x_1) (x_0 + m x_1) & - n (m +
      n) (x_0 + m x_1) x_2 -b x_0 x_2    \\
 &   m^2 (n x_0 - x_1)^2 + b x_2^2 &  m (m + n) (n x_0 - x_1) x_2 - b x_1 x_2  \\
 & &  (m + n)^2 x_2^2+b (x_0^2 + x_2^2)
 \end{array}\!\!\! \right]
\end{eqnarray*} 
is positive semidefinite for all $\x \in \mathbb{P}^2(\mathbb{R})$.
The first two leading principal minors of $\Phi_{b,m,n}(\x \x^T)$,
$$b\, x_1^2+ n^2 (x_0 + m x_1)^2\  \mbox{ and }\ 
b\, m^2 x_1^2 (n x_0 - x_1)^2+ b^2 x_1^2 x_2^2 + b\, n^2 (x_0 + m x_1)^2 x_2^2$$
are  everywhere nonnegative if and only if $b\geq 0$. 
We will study the positivity of $\det \Phi_{b,m,n}(\x \x^T)$ by presenting it as a quadric in $x_2^2$. The determinant 
equals $b^2$ times
\begin{eqnarray} m^2 x_0^2 x_1^2 (n x_0 - x_1)^2  + \label{eq:det8mnpts} \nonumber \\ 
 \left( n^2(1 \!-\! m^2) x_0^4 \!+\! 
   2 m^2 n x_0^3 x_1 \!+\!  (b \!-\! (m \!-\! n)^2) x_0^2 x_1^2 + 2 m^2 n x_0 x_1^3 \!-\! ( b  \!+\! m^2 \!+\!  2 m n\!+\! m^2 n^2) x_1^4\right) x_2^2 +  \nonumber \\
     \left( -(b  + 2 m n+ n^2) x_0^2 - 
   2 m^2 n x_0 x_1 + (b + (m+n)^2 + 
   m^2 n^2) x_1^2 \right) x_2^4\,  \phantom{.+} \nonumber
   \end{eqnarray} 
 and is not constantly 0 for $b\neq 0$.
 Moreover, for $b\ne 0$ it is everywhere nonnegative if and only if
 the constant and the leading term are greater or equal to 0 for all real $x_0,x_1$, and either it has nonnegative middle term
 \begin{equation} \label{eq:middle8}
 n^2(1 \!-\! m^2) x_0^4 + 
   2 m^2 n x_0^3 x_1 + (b \!-\! (m \!-\! n)^2) x_0^2 x_1^2 + 2 m^2 n x_0 x_1^3 -( b  \!+\! m^2 \!+\!  2 m n\!+\! m^2 n^2) x_1^4 \geq 0
   \end{equation} 
   or it has nonpositive discriminant 
 \begin{eqnarray} \label{eq:discriminant8}  (x_0^2 - x_1^2)^2 \times \\
\left( \left(n (m+1) x_0 - m x_1 \right)^2 + (b + 2 m n + 
      m^2 n^2) x_1^2\right) \times \nonumber \\
      \left( \left(n (m-1) x_0 - m x_1\right)^2 + (b + 2 m n + 
      m^2 n^2 ) x_1^2 \right) \leq 0. \nonumber
   \end{eqnarray}
Obviously, the constant term  $m^2 x_0^2 (n x_0 - x_1)^2 x_1^2 $ is nonnegative everywhere.  
The leading term 
\begin{equation} \label{eq:leadtermmn}
-(b  + 2 m n+ n^2) x_0^2 -   2 m^2 n x_0 x_1 + \left( b + (m+n)^2 +   m^2 n^2 \right) x_1^2  \end{equation}
also needs to be always nonnegative, which in particular implies (since $b\geq 0$) that  $m n\leq 0$ and its discriminant 
$4 (b + (m+n)^2 ) (b + 2 m n + n^2 + m^2 n^2)$ is nonpositive.
It turns out that the assumption 
 $$b + 2 m n + n^2 + m^2 n^2 +m^2 \leq m^2$$ 
is sufficient. Indeed, unless $b,m,n$ are all 0,  the discriminant of~\eqref{eq:leadtermmn} is less or equal to 0 and the term at $x_1^2$ is strictly positive, therefore
the leading term of $\det \Phi_{b,m,n}(\x \x^T)$ is  always nonnegative.

Before we proceed with our analysis of positivity, we make some observations about the symmetry of $m$ and $n$.
Recall the dual of a positive map defined in~\eqref{eq:dualhat} by $\y^T\Phi (\x \x^T)\y= \x^T \hat{\Phi} (\y \y^T)\x$. 
In our case
 $\hat{\Phi}_{b,m,n}(\y \y^T)$ equals
$$\left[\!\!\! \begin{array}{ccc}
n^2 (y_0 - m y_1)^2 + b y_2^2& 
 m n (n y_0 + y_1) (y_0 - m y_1) & - b y_0 y_2 + n (m + n) (-y_0 + m y_1) y_2    \\
 &   b y_0^2 + 
  m^2 (n y_0 + y_1)^2 & - b y_1 y_2 - m(m+n) (n y_0 + y_1) y_2  \\
 & &   (m + n)^2 y_2^2 + b (y_1^2 + y_2^2)
 \end{array}\!\!\! \right].$$ 
 It is easy to check that  $\y^T\Phi_{b,m,n} (\x \x^T)\y$ is invariant under the following changes of parameters
 $$\begin{array}{cccc}
 m & n & x_0 & y_0 \\
 \downarrow &\downarrow &\downarrow &\downarrow \\
 -m & -n& -x_0&-y_0
 \end{array}$$
 and
 $$\begin{array}{cccccccc}
 m & n & x_0 & x_1 & x_2 & y_0 & y_1 & y_2 \\
 \downarrow &\downarrow &\downarrow &\downarrow & \downarrow &\downarrow &\downarrow &\downarrow \\
 n & m& y_1 & y_0 & y_2 & x_1 & x_0 & x_2
 \end{array}.$$
 We remark that the latter parameter change can be written as 
 $$\Phi_{b,n,m}(\x \x^T) = 
 \left[ \!\! \begin{array}{ccc} 0 &1 &0 \\ 1&0&0\\ 0&0 &1 \end{array}  \!\! \right] 
 \hat{\Phi}_{b,m,n}\left( 
 \left[ \!\! \begin{array}{ccc} 0 &1 &0 \\ 1&0&0\\ 0&0 &1 \end{array}  \!\! \right]
 \x \x^T 
 \left[ \!\! \begin{array}{ccc} 0 &1 &0 \\ 1&0&0\\ 0&0 &1 \end{array}  \!\! \right]
 \right)
 \left[ \!\! \begin{array}{ccc} 0 &1 &0 \\ 1&0&0\\ 0&0 &1 \end{array}  \!\! \right].
 $$ 
In other words, for each inequality $I(m,n,b)$ that  arises from $\Phi_{b,m,n}(\x \x^T)\succeq 0$ for all $\x$, we also have an analogous inequality $I(n,m,b)$ arising from $\hat{\Phi}_{b,m,n}(\y \y^T) \succeq 0$ for all $\y$. 
For example, we get
\begin{eqnarray*}
b\geq 0,\\
b + 2 m n + n^2 + m^2 n^2 +m^2 \leq m^2,\\
b + 2 m n + n^2 + m^2 n^2 +m^2 \leq n^2,
\end{eqnarray*}
which define region $\mathcal{A}$ shown in Figure~\ref{pictAsymmetry},
$$\mathcal{A}=\left\{ (m,n,b)\in \mathbb{R}^3\colon 0\leq b\leq \min \{-(2 m n + n^2 + m^2 n^2),-(2 m n +  m^2 n^2 +m^2) \}\right\}.$$
\begin{figure}
\begin{center}
\includegraphics[width=5.5cm,height=4cm]{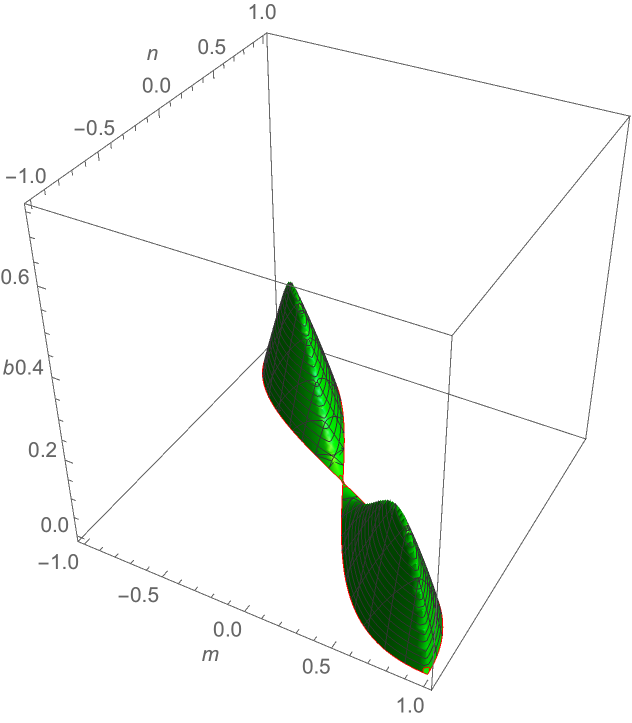}
\includegraphics[width=6.5cm]{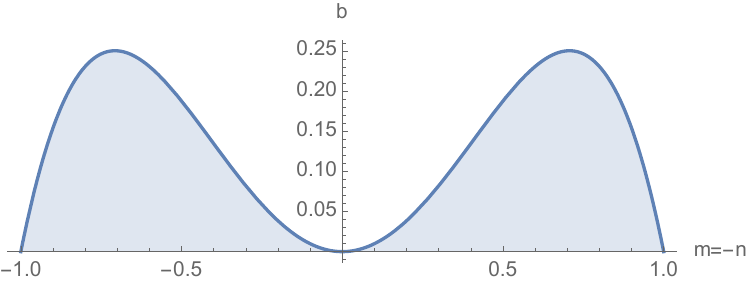}
\end{center}
\caption{ Region $\mathcal{A}$ and its section along $m=-n$.}
\label{pictAsymmetry} 
\end{figure}
An easy calculation shows that the graphs of $f_m(m,n)=-(2 m n +  m^2 n^2 +m^2)$ and $f_n(m,n)=-(2 m n +  m^2 n^2 +n^2)$ defined for $m n\leq 0$ meet along $m=-n$ and that $\min \{ f_m, f_n  \}$ reaches its maximal value $\frac{1}{4}$ at $m=-n=\pm \frac{1}{\sqrt{2}}$. Figure~\ref{pictmnsymmetry} shows the symmetries in $m$ and $n$ when $b=0$. 

Our aim is to prove that  $\Phi_{b,m,n}$ is positive for all $ (m,n,b) \in \mathcal{A}$. Since the cone of positive maps is a convex set and $\Phi_{b,m,n}(\x \x^T)$ is linear in $b$, it is enough to prove that $\Phi_{b,m,n}$ 
is positive for $b=0$ and $b= \min \{-(2 m n + n^2 + m^2 n^2),-(2 m n +  m^2 n^2 +m^2) \}$. Note that $\Phi_{0,m,n}$  is a congruence map, so it is completely positive by definition. 
Proof of positivity of $\Phi_{b,m,n}$ will be finished when we show that 
either \eqref{eq:middle8} or \eqref{eq:discriminant8} holds
for each $x_0,x_1\in\mathbb{R}$ and 
$$b= \min \{-(2 m n + n^2 + m^2 n^2),-(2 m n +  m^2 n^2 +m^2) \}.$$
By the above discussed symmetries of $m$ and $n$ shown in Figures~\ref{pictAsymmetry} and~\ref{pictmnsymmetry}, 
\begin{figure}
\begin{center}
\includegraphics[width=5.5cm]{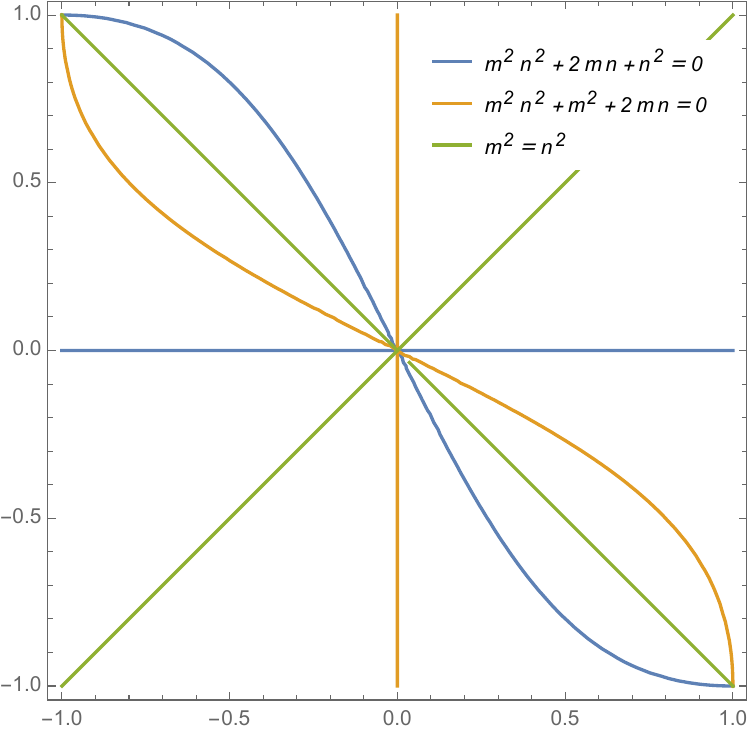}
\includegraphics[width=5.5cm]{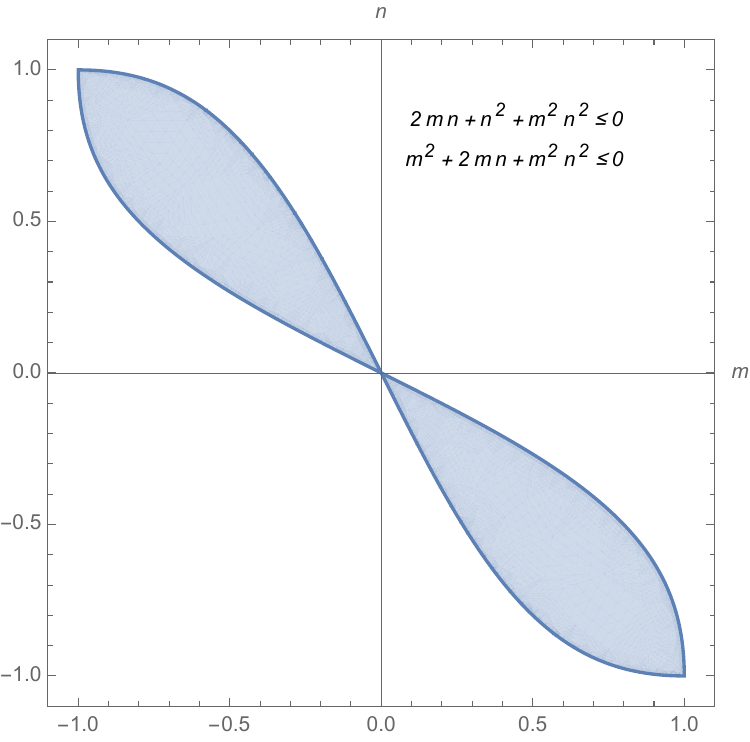}
\end{center}
\caption{Symmetry over the lines $|m|=|n|$.}
\label{pictmnsymmetry} 
\end{figure}
it is enough  to consider 
$$b=-(2 m n +  m^2 n^2 +m^2) \mbox{ and } |m|\geq |n|.$$
Note that in particular this implies $ |m|\leq 1$: if $m > 0$, then $2n+(n^2+1)m\leq 0$ and
$(n^2+1) |n| \leq (n^2+1)m\leq 2 |n|$, which shows that 
$$|n| \leq 1 \mbox{ and } m \leq \frac{-2 n}{n^2 +1}\leq 1.$$
Analogous argument works if $m \leq 0$.  Additionally we assume that $m\neq 0$ and $n \neq 0$, since the contrary would imply $m=n=b=0$.
Under these assumptions we cover $\mathbb{R}^2$ with lines through the origin and show that 
\begin{enumerate}
\item[(i)] inequality \eqref{eq:discriminant8}  holds for points on $x_1 = \frac{(k + m) n}{2 m}x_0$ with $|k|\geq 1$,
\item[(ii)] inequality  \eqref{eq:middle8}  holds for points on  $x_1 = \frac{(k + m) n}{2 m}x_0$ with $|k|\leq 1$,
\item[(iii)] the middle term  \eqref{eq:middle8} equals 0 on the line $x_0=0$.
\end{enumerate}
If we substitute $b=-(2 m n +  m^2 n^2 +m^2) $ and $x_1 = \frac{(k + m) n}{2 m}x_0$  into  \eqref{eq:discriminant8}, we get
$$(1 - k^2) (1 - m^2) n^4 x_0^4 (x_0^2-x_1^2)^2 \leq 0,$$ which proves (i). By the same substitution \eqref{eq:middle8} becomes  
$$\frac{n^2 x_0^4}{4 m^2} \left( 4 m^2 - 2 k^2 m^2 - 2 m^4 - k^2 n^2 - 2 k m n^2 + k^3 m n^2 - 
 m^2 n^2 + 2 k^2 m^2 n^2 + k m^3 n^2 \right).$$
 We claim that this is nonnegative for $|k| \leq 1$, which is equivalent to 
 $$2 m^2 \left(2-k^2-m^2\right) \geq n^2 (k+m)^2(1-k m).$$
 For $|k|=|m|=1$ the inequality clearly holds, else it is equivalent to 
 $$\frac{ m^2}{n^2}\geq \frac{(k+m)^2(1-k m)}{2 \left(2-k^2-m^2\right)}.$$
 The last inequality follows from 
 $$1 \geq \frac{(k+m)^2(1-k m)}{2 \left(2-k^2-m^2\right)},$$
 which is the same as 
 $2 (1-k^2)(1-m^2)+(1-k m)(2-k^2-m^2)\geq 0$ that is clearly true.
 Claim (iii) follows directly from substituting $b=-(2 m n +  m^2 n^2 +m^2) $ and $x_0=0$ in~\eqref{eq:middle8}.
 This concludes the proof of positivity of $\Phi_{b,m,n}$.
 
Finally we will show that the positive maps $\Phi_{b,m,n}$ are not completely positive for all $(m,n,b)\in \mathcal{A}$ with  $b\neq 0$, by proving that the associated nonnegative biquadratic forms $p_{b,m,n}$ have only finitely many (however at least 8) zeros.  It is enough to check that $p_{b,m,n}$ has finitely many real zeros on the boundary of $\mathcal{A}$, since biquadratic forms in the interior of $\mathcal{A}$ are convex combinations of the ones on the boundary.  By symmetry we choose $b=-(2 m n+n^2+m^2 n^2)$ and $|n| \geq |m|$. 
We remark that we were not able to solve the system
$$\frac{\partial \det \Phi_{b,m,n}(\x \x^T)}{\partial x_i}=0,\ i=0,1,2 $$
and apply the same methods as in Section~\ref{sec:9zero}.
Instead we will prove that 
 $\det \Phi_{b,m,n}(\x \x^T) $ and $\det \hat{\Phi}_{b,m,n}(\y \y^T) $ have only finitely many zeros in $ \mathbb{P}^2(\mathbb{R})$, which will  imply that  $p_{b,m,n}$ has only finitely many zeros in $\mathbb{P}^2(\mathbb{R})\times \mathbb{P}^2(\mathbb{R})$ . By~\cite[Theorem 3.5]{CLR}  a nonnegative ternary form has infinitely many real zeros if and only if it is divisible by the square of some indefinite form. By contradiction we will assume that $\det \Phi_{b,m,n}(\x \x^T)$ is divisible by the square of some indefinite form. We will leave an analogous proof that  $\det \hat{\Phi}_{b,m,n}(\y \y^T)$ is not divisible by the square of some indefinite form to the reader.
For  $b=-(2 m n+n^2+m^2 n^2)\neq 0$  we get that
$\det \Phi_{b,m,n}(\x \x^T) $ is $b^2$ times
\begin{eqnarray} m^2 x_0^2 x_1^2 (n x_0 - x_1)^2  + \label{eq:det8bmnpts} \nonumber \\ 
 \left(  n^2 (1- m^2) x_0^4   + 2 m^2 n x_0^3 x_1 - (m^2 n^2+m^2+2 n^2) x_0^2 x_1^2  + 2 m^2 n x_0 x_1^3 +(n^2- m^2) x_1^4  \right) x_2^2 +  \nonumber \\
      m^2 (n x_0 - x_1)^2   x_2^4\, , \phantom{+} \nonumber
   \end{eqnarray} 
which we view as quartic in $x_2$ with coefficients in $\mathbb{R}[x_0,x_1]$. This nonnegative quartic then needs to have one of the following forms
\begin{enumerate}
\item[(I)] it is a nontrivial product of polynomials of degree 4 and 0 in $x_2$,
\item[(II)] it equals a product of two different nonnegative quadrics in $x_2$,
\item[(III)] it is a square of a quadric in $x_2$. 
\end{enumerate}
The coefficient at $x_2^4$ in the first factor in (I) needs to be everywhere nonnegative and of degree strictly less than 2, therefore it can only be a constant. This implies that $(n x_0 - x_1)^2$ divides $\det \Phi_{b,m,n}(\x \x^T)$. If we substitute $x_1=n x_0$ into $\det \Phi_{b,m,n}(\x \x^T)$ we get $b^2 n^2 \left( n^2-1\right)^2  x_0^4 x_2^2$, which is not constantly 0 unless $n=0$ or $n=\pm 1$. But then also $b=0$, which is a contradiction.
The product in (II) equals 
\begin{equation} \label{eq:prodII}
  \left(  m^2 (n x_0 - x_1)^2   x_2^2 +P(x_0,x_1) x_2+Q(x_0,x_1)  \right)\left( x_2^2 +R(x_0,x_1) x_2+S(x_0,x_1)  \right)
  \end{equation}
for some $P,Q,R,S \in \mathbb{R}[x_0,x_1]$, since the leading terms at $x_2$ in both factors need to be everywhere nonnegative. From the coefficients at $x_2^3$ and $x_2$ being 0, we get
$P(x_0,x_1)= -m^2 (n x_0 - x_1)^2 R(x_0, x_1)$ and
$$R(x_0,x_1) \left( Q(x_0,x_1)-m^2 (n x_0-x_1)^2  S(x_0,x_1) \right)=0.$$
Then either $Q(x_0,x_1)=m^2 (n x_0-x_1)^2  S(x_0,x_1)$ or $R(x_0,x_1)=0$. 
If $Q(x_0,x_1)=m^2 (n x_0-x_1)^2  S(x_0,x_1)$, the first factor in~\eqref{eq:prodII} and subsequently $\det \Phi_{b,m,n}(\x \x^T)$  are divisible by $(n x_0-x_1)^2$, which we already considered in (I).
If $R(x_0,x_1)=0$ then also $P(x_0,x_1)=0$. Since in this case $Q(x_0,x_1)$ is not divisible by $(n x_0 - x_1)^2$ (else we get (I) again) and 
 $Q(x_0,x_1) S(x_0,x_1)\neq 0$, neither of the factors in~\eqref{eq:prodII} can be a square of some indefinite form.  
 In the remaining case (III) we assume that  
$$\det \Phi_{b,m,n}(\x \x^T)= b^2
 \left(  m (n x_0 - x_1)   x_2^2 +T(x_0,x_1) x_2\pm m x_0  (n x_0 - x_1) x_1 \right)^2$$
 for some $T(x_0,x_1)$.
 Since $\det \Phi_{b,m,n}(\x \x^T)$ contains no odd powers of $x_2$, the quadric $T(x_0,x_1)$ is constantly 0 and $(n x_0 - x_1)^2$ divides $\det \Phi_{b,m,n}(\x \x^T)$, which as before leads to case (I).

We conclude the proof of Theorem~\ref{thm8zeros} by determining which of  $\Phi_{b,m,n}$ are extremal positive maps. For fixed nonzero $m$ and $n$  we found all nonnegative biquadratic forms  $p_{b,m,n}$ with the set of 8 zeros prescribed at the beginning of Section~\ref{sec:8zero}, where $b$ is an arbitrary number between 0 and $\min \{-(2 m n + n^2 + m^2 n^2),-(2 m n +  m^2 n^2 +m^2) \}$. Since $p_{b,m,n}$ is linear in $b$, the extremal  $\Phi_{b,m,n}$ are exactly the ones with $(m,n,b)$ on the boundary of $\mathcal{A}$.
 
\end{proof}
\medskip

\section{Nonnegative biquadratic forms with Choi set of 7 zeros}
\label{sec:7zero}

Nonnegative biquadratic forms with Choi set of zeros
$$\begin{array}{l}\vspace{1mm}
 \left\{  (1,1,1;1,1,1), (1,1,-1;1,1,-1), (1,-1,1;1,-1,1), (-1,1,1;-1,1,1), \right. \\
\  \, \left.  (1,0,0; 0,1,0), (0,1,0;0,0,1), (0,0,1;1,0,0)\right\},
\end{array}$$ 
had previously been determined in~\cite{Q}. Using the canonical form from Section~\ref{SecCan}, we are able to present the corresponding positive maps and the area of positivity in a more symmetric way. 

\begin{figure}
\begin{center}
\includegraphics[width=8cm]{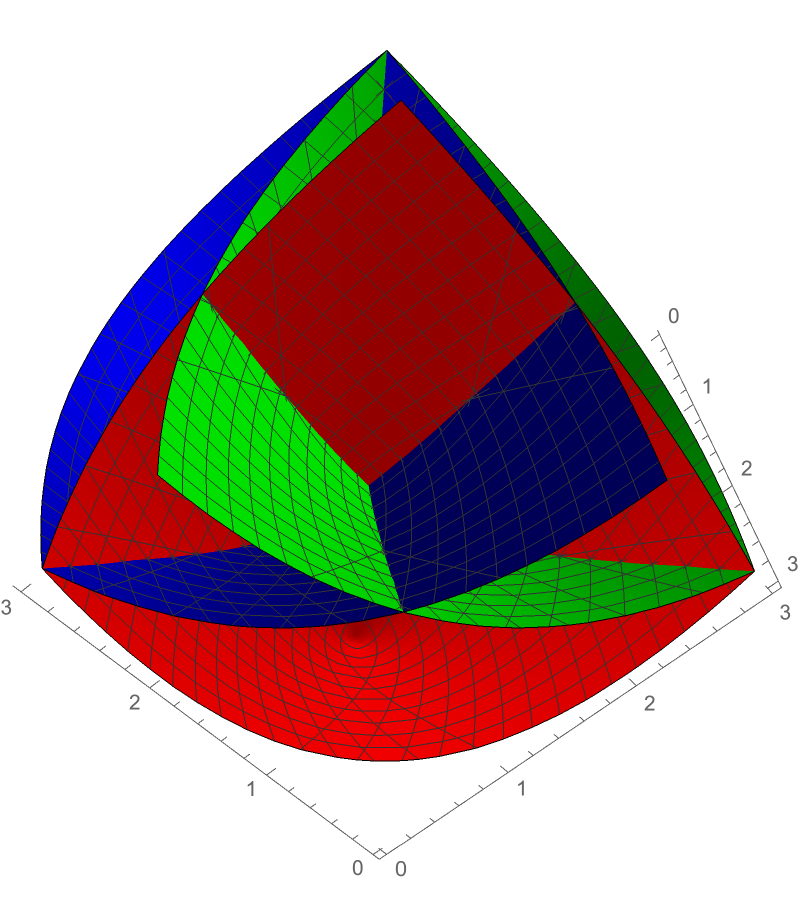}
\end{center}
\caption{The  $2 \times 2$ principal minors imply semidefiniteness}
\label{figChoi7}
\end{figure}

In Sections \ref{sec:10zero}, \ref{sec:9zero} and \ref{sec:8zero} we found that for 
nonnegative biquadratic forms $p_0$ with the Choi set of 7 zeros,
the associated $\Phi_0(\x \x^T)$ are of the form 
$$\left[\!\! \begin{array}{ccc}
a_{00} x_0^2 +  a_{11} x_1^2 & -\frac{1}{2} ( a_{11} +a_{00} + d_{11} - f_{22}) x_0 x_1 & -\frac{1}{2} ( a_{11} +a_{00}  - d_{11} + f_{22}) x_0 x_2\\
& d_{11} x_1^2 + a_{11} x_2^2 &  -\frac{1}{2} ( a_{11} -a_{00}  + d_{11} + f_{22}) x_1 x_2 \\
 &  & a_{11} x_0^2 + f_{22} x_2^2 
  \end{array}\!\! \right].$$
We omit detailed analysis when  $\Phi_0(\x \x^T)$ is positive semidefinite for all $\x \in \mathbb{P}^2(\mathbb{R})$, since it is similar to the methods used in~\cite[Section 6.4]{Q} and in our previous sections. We view the principal minors as linear, quadratic and cubic in $x_0^2,x_1^2,x_2^2$, and then show that the nonnegativity of the $1\times 1$ and $2 \times 2$ principal minors implies the nonnegativity of the determinant. 
When $a_{11}= 0$, the above determinant is 0 and the associated biquadratic form  $p_0$ is SOS under the following conditions
\begin{eqnarray*}
a_{00}\geq 0,\ d_{11}\geq 0,\ f_{22}\geq 0,\\
-a_{00}^2-d_{11}^2- f_{22}^2 +2 a_{00}  d_{11} +2 a_{00} f_{22} + 2  d_{11}  f_{22} \geq 0.
\end{eqnarray*}
When $a_{11}\neq  0$, we can without loss of generality take $a_{11}=1$.
Then   $\Phi_0(\x \x^T)$ is positive semidefinite for all $\x \in \mathbb{P}^2(\mathbb{R})$ 
under the following conditions
\begin{eqnarray} \label{ChoiCond}
a_{00}\geq 0,\ d_{11}\geq 0,\ f_{22}\geq 0,\nonumber \\
 4 d_{11} f_{22}-(1-a_{00}  + d_{11} + f_{22})^2 \geq 0, \\
 4 a_{00} d_{11}-(1+a_{00}  + d_{11} - f_{22})^2 \geq 0, \nonumber \\
 4 a_{00} f_{22}-(1+a_{00}  - d_{11} + f_{22})^2 \geq 0.\nonumber
\end{eqnarray}
The area cut out by the above three quadratic surfaces is shown on Figure~\ref{figChoi7}.
This way we also obtain a big family of positive maps that are not completely positive. The special case $a_{00}=d_{11}=f_{22}$ has already been studied in~\cite[Section 5.3.2]{Q}, and  for $a_{11}=a_{00}=d_{11}=f_{22}$
we get the  example of Choi 
$$\left[\!\! \begin{array}{ccc}
 x_0^2 +   x_1^2 & -x_0 x_1 & - x_0 x_2\\
&  x_1^2 +x_2^2 &   -x_1 x_2 \\
 &  &x_0^2 + x_2^2 
  \end{array}\!\! \right]$$ whose determinant is  
  $x_0^4 x_2^2 + x_1^2 x_2^4 + x_0^2 x_1^4 - 3 x_0^2 x_1^2 x_2^2$. 

\section{Extremal positive maps on $M_3(\mathbb{C})$}
\label{SecExtension}

The aim of this section is the following. For each of the extremal positive linear maps on $\Sym_3$ considered in  Sections~\ref{sec:10zero},~\ref{sec:9zero} and~\ref{sec:8zero}, we will present an extension to $M_ 3 (\mathbb{C})$ that will be extremal in the cone of positive maps on $M_ 3 (\mathbb{C})$. 
A positive linear map $\Phi \colon \Sym_n \rightarrow \Sym_m$ can be extended to a positive linear map $\Phi \colon M_n(\mathbb{C}) \rightarrow M_m(\mathbb{C})$ in various ways. 
Let $\{ E_{ij};0\le i,j\le n-1 \}$ be the standard basis for $M_n(\mathbb{C})$; then $E_{ij}+E_{ji}$ and $E_{ij}-E_{ji}$ define symmetric and skew-symmetric matrices respectively, for $0\leq i <j \leq n-1$. Since $\Phi$ is linear over $\mathbb{C}$,
the equality 
$$Z=\left[ z_{ij} \right]_{i,j=0}^{n-1}=\frac{1}{2} \left(Z+Z^{T} \right)+\frac{1}{2} \left(Z-Z^{T} \right)$$ implies
$$\Phi(Z)=
\sum_{i=0}^{n-1} z_{ii} \, \Phi (E_{ii})+\sum_{0\leq i <j \leq n-1}\left(\frac{z_{ij}+z_{ji}}{2} \,  \Phi (E_{ij}+E_{ji})+\frac{z_{ij}-z_{ji}}{2}\,  \Phi (E_{ij}-E_{ji})\right).$$

In~\cite{KMcCSZ} the authors extended $\Phi$ from $\Sym_n$ to $M_n(\mathbb{R})$ by defining it as the zero map on the skew-symmetric matrices, and then proved that its complexification yields a positive linear map on $M_n(\mathbb{C})$. However, an extremal positive map on $\Sym_n$ will most likely not extend to an extremal positive map on $M_ n (\mathbb{C})$. For example, the following extension of the Choi map 
$$\begin{array}{ccc}
Z
 &\mapsto &
\left[\!\! \begin{array}{ccc}
 z_{00} +    z_{11} & -\frac{1}{2} \left( z_{01}+z_{10} \right) &-\frac{1}{2} \left( z_{02}+z_{20} \right) \\
-\frac{1}{2} \left( z_{01}+z_{10} \right)&  z_{11} +  z_{22} &   -\frac{1}{2} \left( z_{12}+z_{21} \right) \\
-\frac{1}{2} \left( z_{02}+z_{20} \right) &  -\frac{1}{2} \left( z_{12}+z_{21} \right) &  z_{00} + z_{22} 
  \end{array}\!\! \right]
   \end{array}$$
is not extremal, since it equals $\frac{1}{2} \left( \Psi+\Psi\circ T\right)$ where
\begin{equation} \label{ChoiHa}
\begin{array}{ccc}
\Psi \colon Z &\mapsto &
\left[\!\! \begin{array}{ccc}
 z_{00} +    z_{11} & -z_{01} & - z_{02} \\
-z_{10} &  z_{11} +  z_{22} &   -z_{12} \\
-z_{20} & -z_{21}  &  z_{00} + z_{22} 
  \end{array}\!\! \right]
   \end{array}\end{equation}
and $T$ is the transpose map.
On the other hand, K.-C. Ha~\cite{H1} proved that $\Psi\colon M_3(\mathbb{C}) \to M_3(\mathbb{C})$ is an extremal positive map.

In this section we will use the same method as in the construction of (extremal) positive maps on $\Sym_3$, which were obtained from the associated biquadratic forms with the prescribed sets of zeros.
For a linear map $\Phi \colon M_ 3 (\mathbb{C}) \to M_ 3 (\mathbb{C})$  we define
$$p_{\Phi}(\x,\y)=\y^{\ast}\Phi \left( \x \x^{\ast} \right) \y,$$
where $ \x = (x_0,x_1,x_2)$ and  $\y =(y_0,y_1,y_2) $ are in $\mathbb{P}^{2}(\mathbb{C})$ and $\ast$ denotes conjugate transpose. 
The map $\Phi$ is positive if and only if
$p_{\Phi}$ is a nonnegative polynomial in variables 
$$\frac{1}{2} \left(x_k+\overline{x_k} \right), \frac{1}{2i} \left(x_k-\overline{x_k} \right)\ \mbox{ and }
\ \frac{1}{2} \left(y_k+\overline{y_k} \right), \frac{1}{2i} \left(y_k-\overline{y_k} \right)$$
in  $\mathbb{P}^{5}(\mathbb{R}) \times \mathbb{P}^{5}(\mathbb{R})$. 
As in Section~\ref{SecCan}, a positive map $\Phi$ comes in pair with another positive map $\hat{\Phi}$, which is defined by the duality
$$ \y^{\ast}\Phi (\x \x^{\ast})\y= \x^{\ast} \hat{\Phi} (\y \y^{\ast})\x. $$
Next, write $\Phi$ in coordinates  $\Phi\left(\left[ z_{ij} \right]_{i,j=0}^2\right)=\left[ \sum_{k,l=0}^2 a_{ijkl} z_{kl} \right]_{i,j=0}^2$, where $a_{ijkl}\in \mathbb{C}$. Then $\hat{\Phi}\left(\left[ z_{ij} \right]_{i,j=0}^2\right)=\left[ \sum_{k,l=0}^2 a_{lkji} z_{kl} \right]_{i,j=0}^2$.
Since positive maps send Hermitian matrices into Hermitian matrices and  since $i (E_{kl}-E_{lk})$ are Hermitian, we obtain the following linear system of equations in $\operatorname{Re}(a_{ijkl})$ and  $\operatorname{Im}(a_{ijkl})$,
$$\Phi(E_{kl})=\Phi(E_{lk})^{\ast},\quad  \hat{\Phi}(E_{kl})=\hat{\Phi}(E_{lk})^{\ast},$$
whose solution is defined by $a_{ijkl}=\overline{a_{jilk}}$ for $0\le i,j,k,l\le 2$. In particular, $a_{kkll}\in \mathbb{R}$ for all $k$ and $l$.
Each zero $(\x, \y)\in \mathbb{P}^{2}(\mathbb{C}) \times \mathbb{P}^{2}(\mathbb{C})$ of $p_{\Phi}(\x,\y)$ imposes additional linear relations among $a_{ijkl}$, since by positivity the following holds
\begin{equation} \label{eq:linzeros}
\Phi \left( \x \x^{\ast} \right) \y=(0,0,0)\  \mbox{ and }\  \hat{\Phi} \left( \y \y^{\ast} \right) \x=(0,0,0).\end{equation}
Therefore, by prescribing a sufficiently large set of zeros for $p_{\Phi}(\x,\y)$ and solving the above linear equations, we obtain candidates for positive maps $\Phi$.

\begin{prop}  \label{propComplex}
The following linear maps on  $M_3(\mathbb{C})$ 
are extremal positive maps. 

\noindent
$\Ia .$ \vspace{-5mm}
$$
\begin{array}{c}
Z \\
\begin{tikzcd}
  \arrow[mapsto]{d}
  \\
\phantom{a} 
\end{tikzcd}\\
\left[
\begin{array}{ccc}(t^2-1)^2 z_{00}+z_{11}+t^4 z_{22}&-(t^4-t^2+1)z_{01}&-(t^4-t^2+1) z_{02} \\
-(t^4-t^2+1) z_{10} &  t^4 z_{00}+(t^2-1)^2 z_{11}+z_{22} &-(t^4-t^2+1) z_{12}\\
-(t^4-t^2+1) z_{20} & -(t^4-t^2+1) z_{21} & z_{00}+t^4 z_{11}+ (t^2-1)^2 z_{22}
\end{array}
\right]
\end{array}$$  
for $ t\in \mathbb{R} \backslash \{\pm1\}$;\\

\noindent
$\Ib.$ \vspace{-8mm}
$$
\begin{array}{c}
Z \\
\begin{tikzcd}
  \arrow[mapsto]{d}
  \\
\phantom{a} 
\end{tikzcd}\\
\!\!\!\! \left[\!\!\! \begin{array}{ccc}
 \substack{ p^2(p q-1)^2} z_{00} + \substack{q(2 p - q) }z_{11} \!\!\!\!\!\!\! & \!\!\!\!\!\!\! \substack{ -p q (1 -  q^2 + p^2 q^2) }z_{01}  
 \!\!\!\!\!\!\! & \!\!\!\!\!\!\! 
\substack{ (p q-1)(p^2 + p q - p^3 q - q^2 + p^2 q^2)} z_{02}   \\
\!\!\!\!\!\!\!\!\!\!\!\!\!\! \substack{ -p q (1 -  q^2 + p^2 q^2) }z_{10}  
 \!\!\!\!\!\!\!  & \!\!\!\!\!\!\!  \substack{ p^2 q^3 (2 p - q)}  z_{00} \substack{+q^2 ( p q-1)^2} z_{11} \substack{+q(2 p-q) } z_{22}    \!\!\!\!\!\!\! & \!\!\!\!\!\!\!  \substack{ -p q (1 -  q^2 + p^2 q^2)} z_{12}  \\
\substack{ (p q-1)(p^2 + p q - p^3 q - q^2 + p^2 q^2)} z_{20} \!\!\!\!\!\!\!  &\!\!\!\!\!\!\!  \substack{ -p q (1 -  q^2 + p^2 q^2)} z_{21} \!\!\!\!\!\!\!\! & \!\!\!\!\!\!\!\!\!\!\!\!\!  \substack{q (2 p - q) (1 - p^2 q^2) } z_{00} \substack{+p^2 q^3 (2 p - q) } z_{11}  \substack{+p^2(p q-1)^2 }z_{22}
 \end{array}\!\!\! \right]
\end{array}$$  
 for 
$$\left\{ (p,q)\in [0,1/\sqrt{2}] \times [0,\sqrt{2}] \colon 2 p - q \ge 0,\ (p^2 - 1)^2 q^2 - p^2\geq 0 \right\}\backslash \left\{(0,0),\left(\frac{1}{\sqrt{2}},\sqrt{2}\right)\right\};$$

\noindent
$\Ic.$
$$\begin{array}{l} \vspace{2mm} 
\hspace{1cm} Z \ \ \  \mapsto \ \  \  b \left[\!\!\! \begin{array}{ccc}
 z_{11} & 0 & -z_{02} \\
 0 & z_{22} & - z_{12}  \\
 -z_{20} & -z_{21} & z_{00} + z_{22}
 \end{array}\!\!\! \right]+
 c \left[\!\!\! \begin{array}{ccc}
 0 &  z_{01}-z_{10} &0 \\
  z_{10}-z_{01}  & 0 & 0 \\
 0 & 0 & 0
 \end{array}\!\!\! \right]+  \\ 
 \! \left[\!\!\! \begin{array}{ccc}
n^2 (z_{00} + m (z_{01}+z_{10})+m^2 z_{11}) & - m n \left(n z_{00}-z_{01}+mn z_{10}-m z_{11}\right) & - n (m +
      n) (z_{02} + m z_{12})      \\
- m n \left(n z_{00}-z_{10}+mn z_{01}-m z_{11}\right)&   m^2 (n^2 z_{00} - n (z_{01}+z_{10})+z_{11}) &  m (m + n) (n z_{02} - z_{12})    \\
 - n (m +
      n) (z_{20} + m z_{21}) &    m (m + n) (n z_{20} - z_{21})  &  (m + n)^2 z_{22}
 \end{array}\!\!\! \right]
 \end{array}$$
for $(m,n)\in [-1,1]^2\backslash \{(0,0)\}$, such that
$$b=\min \left\{-\! 2 m n \!-\! m^2 n^2 \!-\! n^2,\, -\! 2 m n \!-\!  m^2 n^2 \!-\! m^2 \right\}\ge 0$$ and 
$$c=\max  \left\{  -\frac{1}{2}  m n (1\!+\!m n \!-\!\sqrt{1 \!-\! m^2}),\, -\frac{1}{2}  m n (1\!+\!m n \!-\! \sqrt{1 \!-\! n^2}),   \right\}.$$
\end{prop}

\begin{proof}
\noindent $\Ia .$ \hspace{3mm}
For a given $t \in \mathbb{R} \backslash \{\pm1,0 \}$, assume that $p_{\Phi_t}(\x, \y)$ is nonnegative with the following set of zeros
$$\begin{array} {l}\left\{  \left(e^{i \varphi_0},e^{i \varphi_1}, e^{i \varphi_2}; e^{i \varphi_0},e^{i \varphi_1}, e^{i \varphi_2} \right), \right. \\ \vspace{1mm}
\  \, \left. \left(1, t\, e^{i \varphi}, 0; t\, e^{-i \varphi}, 1, 0 \right), 
\left(0, 1, t\, e^{i \varphi}; 0, t\, e^{-i \varphi}, 1 \right),
\left(t\, e^{i \varphi}, 0, 1; 1, 0, t\, e^{-i \varphi} \right)
\right\},
\end{array}$$
for all $\varphi, \varphi_0, \varphi_1, \varphi_2 \in [0,2 \pi)$.
It is straightforward to check that the  linear system of equations, obtained from evaluating~\eqref{eq:linzeros} at these zeros, has a unique solution up to a scalar,
which defines the map $\Phi_t$ as stated in the proposition.
Positivity of $\Phi_t$ follows from calculating the trace of $\Phi_t(\x \x^{\ast})$
$$
2 (1 - t^2 + t^4) (|x_0|^2  + |x_1|^2 + |x_2|^2 ),$$
the sum of the main $2\times 2$ minors of $\Phi_t(\x \x^{\ast})$
 $$
(1 - t^2 + t^4)^2 (|x_0|^2  + |x_1|^2 + |x_2|^2)^2$$
and the determinant of $\Phi_t(\x \x^{\ast})$
$$\begin{array}{c}
   (1-t^2)^2\times  \hspace{9cm}\\
    \left(   t^4 (|x_0|^6  + |x_1|^6 + |x_2|^6) +
    (t^8 - 2 t^2) (|x_0|^4 |x_1|^2 + 
   |x_0|^2 |x_2|^4    + 
   |x_1|^4  |x_2|^2 )+ \right. \\
  \left.  (1 - 2 t^6) (|x_0|^2  |x_1|^4 + 
   |x_0|^4  |x_2|^2 + 
   |x_1|^2 |x_2|^4)
   -3 (1 - 2 t^2 + t^4 - 2 t^6 + t^8) |x_0|^2|x_1|^2|x_2|^2  \right)
   \end{array}$$
which can be considered as polynomials in $|x_0|^2, |x_1|^2, |x_2|^2$. As such they are equal to the corresponding polynomials in Section~\ref{sec:10zero} and are therefore all nonnegative.

Finally, $\Phi_t$ is extremal since it is up to a scalar the only positive linear map with the prescribed set of zeros. 
For $t=0$ we get the Choi map that is extremal by~\cite{H1}.

\noindent $\Ib .$ \hspace{3mm}
For given nonzero $p,q \in \mathbb{R}$ with $pq\ne 1$, we find all nonnegative polynomials $p_{\Phi_{p,q}}(\x, \y)$ with the following set of zeros
$$\begin{array} {l}\left\{  \left(e^{i \varphi_0},e^{i \varphi_1}, e^{i \varphi_2}; e^{i \varphi_0},e^{i \varphi_1}, e^{i \varphi_2} \right), \right. \\ \vspace{1mm}
\  \, \left. \left(1, p\, e^{i \varphi}, 0; q\, e^{-i \varphi}, 1, 0 \right), 
\left(0, 1, q\, e^{i \varphi} ; 0, p\, e^{-i \varphi}, 1 \right),
\left(0, 0, 1; 1, 0, 0 \right)
\right\},
\end{array}$$
for all $\varphi, \varphi_0, \varphi_1, \varphi_2 \in [0,2 \pi)$. 
The linear system of equations, obtained from evaluating~\eqref{eq:linzeros} at these zeros,  has (up to a scalar) a unique solution that defines $\Phi_{p,q}$ as stated in the proposition.

We will check that $\Phi_{p,q}\colon M_3(\mathbb{C}) \to M_3(\mathbb{C}) $ is positive for the same pairs $(p,q)$ as its restriction  $\Phi_{p,q}\colon \Sym_3 \to \Sym_3$ considered in Section~\ref{sec:9zero}. The first and second 
principal minors of $\Phi_{p,q}(\x \x^{\ast})$ are 
\begin{eqnarray*}
 p^2(p q-1)^2 |x_0|^2 + q(2 p - q) |x_1|^2\  \mbox{ and}\\
 (2 p - q) q  \left( q^2 (1 - p q)^2 (p^2 |x_0|^2  - |x_1|^2)^2  + p^2 (1 - p q)^2 |x_0|^2  |x_2|^2 +
     q (2 p - q) |x_1|^2 |x_2|^2 \right), \hspace{1cm}
 \end{eqnarray*}
and its determinant equals 
$q^2( p q-1)^2 (2 p - q)^2$ times 
\begin{eqnarray*} p^4 q^2 (1 - p^2 q^2) |x_0|^6 + p^2 q^4  |x_1|^6+\\
 p^2 q^2 (-2 + 2 p^2 q^2 + p^4 q^2) |x_0|^4 |x_1|^2+ p^2 (1 - 2 p^2 q^2 + q^4 - 2 p^2 q^4 + p^4 q^4) |x_0|^4 |x_2|^2+\\ 
q^2 (1 -  p^2 q^2 - 2 p^4 q^2) |x_0|^2 |x_1|^4- ( p^2 + 2 q^2 - 6 p^2 q^2 + p^2 q^4 - 
  2 p^4 q^4 + p^6 q^4) |x_0|^2 |x_1|^2 |x_2|^2+ \\
  +(q^2 -p^2+ p^4 q^2 -  2 p^2 q^2) 
|x_0|^2 |x_2|^4 - 2 p^2 q^2 |x_1|^4 |x_2|^2+p^2 |x_1|^2 |x_2|^4.
\end{eqnarray*}
The three principal minors are nonnegative under the same conditions as in the real case considered in the proof of Theorem~\ref{thm9zeros}. For each $(p,q)$ in the region shown on Figure~\ref{pictRegion9pts}, except for $(p,q)\in \{(0,0),(\frac{1}{\sqrt{2}},\sqrt{2})\}$, $\Phi_{p,q}$ is up to a scalar the only positive linear map associated to a polynomial with the prescribed set of complex zeros, therefore $\Phi_{p,q}$  is extremal.

\noindent $\Ic .$ \hspace{3mm}
We start by observing that
$$2mn+m^2n^2+\max\{m^2,n^2\}=\big(1+mn-\sqrt{1-\max\{m^2,n^2\}}\big)\big(1+mn+\sqrt{1-\max\{m^2,n^2\}}\big),$$
with the second factor on the right hand side being strictly positive, unless $m=\pm 1,n=\mp 1$, in which case  the first factor is also zero. Therefore $b=0$ if and only if $c=0$. Moreover, in this case $\Phi_{m,n}$  is a congruence map, and thus positive and extremal. In the following we will therefore assume that $2mn+m^2n^2+\max\{m^2,n^2\}\ne 0$, and in particular $(m,n)\ne (\pm 1,\mp 1)$.

For fixed nonzero $m, n \in  \mathbb{R}$ with $mn\ne -1$ and all $\varphi \in [0,2 \pi)$ we prescribe the following zeros of $p_{\Phi_{m,n}}(\x, \y)$,
$$\begin{array} {l}\left\{  \left(1,1, e^{i \varphi}; 1,1, e^{i \varphi} \right),\,  \left(1,-1, e^{i \varphi}; 1,-1, e^{i \varphi} \right), \right. \\ \vspace{1mm}
\  \, \left. 
\left(1, 0, 0; m, 1, 0 \right), 
\left(1, n, 0; 0, 1, 0 \right),
\left(0, 1, 0; 0, 0, 1 \right), \left(0, 0, 1; 1, 0, 0 \right)
\right\}.
\end{array}$$
If $m\ne -n$, then the solution of the linear system of equations, obtained from evaluating~\eqref{eq:linzeros} at these zeros,  defines (up to a scalar) $\Phi_{m,n}$ with real parametres $b$ and $c$ as stated in the proposition. It remains to determine for which
 $b, c  \in \mathbb{R}$ the map $\Phi_{m,n}$  is positive and extremal. Note that $\Phi_{m,n}$ restricted to $\Sym_3$ equals to the map in Theorem~\ref{thm8zeros}, which is positive for $0\le b\le \min\{-(2mn+n^2+m^2n^2),-(2mn+m^2+m^2n^2)\}$ and extremal for 
$b= \min \left\{-(2 m n + n^2 + m^2 n^2),-(2 m n +  m^2 n^2 +m^2) \right\}$. The same as in Section~\ref{sec:8zero} (see Figure~\ref{pictmnsymmetry}), it is enough to consider $m$ and $n$ for which 
\begin{equation}1> m\geq -n >0.
\label{eq:assumptionmn} \end{equation} 
Assume now that $p_{\Phi_{m,n}}$ vanishes also at the points
\begin{equation}\label{eq:additional_zeros}
(e^{i\varphi},n\cos\varphi,0;m\sin\varphi,\sin\varphi+ i \sqrt{1-m^2} \cos\varphi,0)
\end{equation}
for all $\varphi \in [0,2\pi)$. The solution of the linear system, obtained from additionally evaluating ~\eqref{eq:linzeros} at these points gives
$$b=-(m^2+2mn+m^2n^2)\quad \mathrm{and}\quad c=-\frac{1}{2}mn(1+mn-\sqrt{1-m^2}),$$
which means that for $1> m\geq -n >0$  the map $\Phi_{m,n}$ is unique (up to a scalar) with $b$ and $c$ as above (note that  $\Phi_{m,n}$ is unique also for $m=-n$). 

We will now show that the obtained map is positive; then it is automatically extremal, as it is equal to the unique positive map associated to the polynomial with the prescribed set of complex zeros. We will verify that $\Phi_{m,n}(\x \x^{\ast})$ is positive semidefinite for all 
$ \x=(r_0 e^{i \varphi_0}, r_1 e^{i \varphi_1},r_2 e^{i \varphi_2})\in \mathbb{P}^{2}(\mathbb{C})$. Its first and second principal minors are
\begin{eqnarray*}n^2r_0^2+2mn^2r_0r_1\cos(\varphi_0-\varphi_1)-m(m+2n)r_1^2&=&
\\n^2\big(r_0+m r_1\cos(\varphi_0-\varphi_1)\big)^2-(m^2+2mn+m^2n^2)^2r_1^2+m^2n^2r_1^2\sin^2(\varphi_0-\varphi_1)&\ge& 0
\end{eqnarray*}
and
\begin{eqnarray*}
-m^2(m^2+2mn+m^2n^2)\big(r_1^2-2nr_0r_1\cos(\varphi_0-\varphi_1)+n^2r_0^2\cos^2(\varphi_0-\varphi_1)\big)r_1^2&-&\\
(m^2+2mn+m^2n^2)\big(n^2r_0^2+2mn^2r_0r_1\cos(\varphi_0-\varphi_1)-(m^2+2mn)r_1^2\big)r_2^2&-&\\
(4c^2+4mn(1+mn)c+m^2n^2(m^2+2mn+m^2n^2))r_0^2r_1^2\sin^2(\varphi_0-\varphi_1)&\ge& 0,
\end{eqnarray*}
where we used the fact that $4c^2+4mn(1+mn)c+m^2n^2(m^2+2mn+m^2n^2)=0$ and $-(m^2+2mn)\geq m^2 n^2$.

It remains to show the positivity of  $\det \Phi_{m,n}(\x \x^{\ast})$, which equals
\begin{eqnarray} \label{eq:detcx}
(m^2 + 2 m n + m^2 n^2)^2 F(r_0,r_1,r_2,\varphi_0,\varphi_1)+&& \\
 4   (m^2 + 2 m n + m^2 n^2) r_0^2 r_1^2
  \big( (m n-n^2 + 2 m^2 n^2)r_2^2 + (m n + m^2 n^2) r_0^2 \big)  \sin(\varphi_0 - \varphi_1)^2 \, c\, +&& \nonumber \\
 4 r_0^2 r_1^2 \big((m^2  n^2-n^2 )r_2^2+(m^2 + 2 m n + m^2 n^2) r_0^2 \big)  \sin(
  \varphi_0 - \varphi_1)^2\,  c^2, \quad &&   \hspace{0.3cm}  \nonumber 
   \end{eqnarray}
where
\begin{eqnarray*}F(r_0,r_1,r_2,\varphi_0,\varphi_1)&=&\\
n^2 ( r_0^2-r_2^2) (r_0^2 - r_1^2)r_2^2 + 
 m^2 r_0^2 ((r_2^2 - r_1^2)^2 + n^2 (r_2^4 + 3 r_1^2r_2^2 - r_0^2 (r_2^2 - r_1^2)))  &-& \nonumber \\
    2 m^2 n r_0 r_1 (r_2^2 - r_0^2)  (r_2^2 - r_1^2) \cos(\varphi_0 - \varphi_1) - 
    4 m^2 n^2 r_0^2 r_1^2r_2^2 \cos(\varphi_0 - \varphi_1)^2.&&
\end{eqnarray*}
The determinant is not constantly zero, because $b>0$.

We will first show that
\begin{equation}\label{eq:estimate_constant_term}
F(r_0,r_1,r_2,\varphi_0,\varphi_1)\ge m^2n^2r_1^2(r_0^2+r_2^2)^2\sin^2(\varphi_0-\varphi_1).
\end{equation}
By continuity it suffices to show this inequality for strictly positive $r_0, r_1, r_2$ and $r_0 \neq r_2$. For fixed numbers $r_0,r_2$ (distinct and positive), denote $\varphi=\varphi_0-\varphi_1$, and define the function 
$$f(r_1,\varphi)= \frac{F(r_0,r_1,r_2,\varphi,0)}{r_1^2}-m^2n^2(r_0^2+r_2^2)^2\sin^2\varphi,$$
which equals
\begin{eqnarray*} 
 && \big((m^2-n^2+m^2n^2)r_2^2+n^2(1-m^2)r_0^2\big)\frac{r_0^2r_2^2}{r_1^2}+m^2r_0^2r_1^2 \\
 &-& 2m^2nr_0\left(r_2^2-r_0^2\right)\left(\frac{r_2^2}{r_1}-r_1\right)\cos\varphi  \\
&+& \left(n^2-m^2n^2\sin^2\varphi \right)r_2^4-\left(2m^2+n^2-m^2n^2 \left(1-2\cos^2\varphi \right) \right) r_0^2r_2^2+m^2n^2r_0^4\cos^2\varphi. 
\end{eqnarray*}
To prove \eqref{eq:estimate_constant_term} it suffices to show that $f$ is nonnegative.
Since for each $\varphi$,
$$\lim_{r_1\to 0}f(r_1,\varphi)=\lim_{r_1\to \infty}f(r_1,\varphi)=\infty,$$
it is enough to show that $f$ is nonnegative in the critical points where
$\frac{\partial f}{\partial \varphi}= \frac{\partial f}{\partial r_1}=0$.
This means that 
$$\frac{2m^2n \left(r_2^2-r_0^2\right)}{r_1}\left(r_0(r_2^2-r_1^2)-nr_1(r_2^2-r_0^2)\cos\varphi\right)\sin\varphi $$
and 
$$\frac{2r_0}{r_1^3}\left(m^2r_0r_1^4-((m^2-n^2+m^2n^2)r_2^2+n^2(1-m^2)r_0^2)r_0r_2^2+m^2nr_1(r_2^2-r_0^2)(r_2^2+r_1^2)\cos\varphi\right)$$ both need to be 0.
As $r_0\ne r_2$, the first equation implies either $\sin \varphi =0$ or 
$\cos \varphi =\frac{r_0(r_2^2-r_1^2)}{nr_1(r_2^2-r_0^2)}$. 
In the second case $\frac{\partial f}{\partial r_1}=\frac{2n^2(1-m^2)(r_2^2-r_0^2)r_0^2r_2^2}{r_1^3}=0$ yields a contradiction.
The critical points of $f$ therefore all satisfy $\sin \varphi=0$. However,
\begin{eqnarray*}f(r_1,k\pi)&=&\left((m^2-n^2+m^2n^2)r_2^2+n^2(1-m^2)r_0^2\right)\frac{r_0^2r_2^2}{r_1^2}+m^2r_0^2r_1^2 \pm\\
& &2m^2nr_0(r_2^2-r_0^2)\left(\frac{r_2^2}{r_1}-r_1\right)+n^2r_2^4+(-2m^2-n^2-m^2n^2)r_0^2r_2^2+m^2n^2r_0^4,
\end{eqnarray*}
which is nonnegative, as shown in Section~\ref{sec:8zero}. This completes the proof of inequality \eqref{eq:estimate_constant_term}.

The nonnegativity of $\det \Phi_{m,n}$ clearly follows from \eqref{eq:estimate_constant_term} and the following inequality 
\begin{eqnarray}
r_1^2\sin^2(\varphi_0-\varphi_1)\Big((m^2 + 2 m n + m^2 n^2)^2m^2n^2(r_0^2+r_2^2)^2&+&  \nonumber\\
 4   (m^2 + 2 m n + m^2 n^2) r_0^2 
  \big( (m n-n^2 + 2 m^2 n^2)r_2^2 + (m n + m^2 n^2) r_0^2 \big) \, c\, &+& \nonumber \\
 4 r_0^2 \big((m^2  n^2-n^2 )r_2^2+(m^2 + 2 m n + m^2 n^2) r_0^2 \big)\,  c^2\Big)&\ge&0  \nonumber 
   \end{eqnarray}
which needs to be proved.
Consider the above expression as a quadric in $r_0^2$ and $r_2^2$. 
The coefficient at $r_2^4$  is clearly positive, and the coefficient at $r_0^4$ is zero by the definition of $c$. It remains to show that
$$4(m^2n^2-n^2)c^2+4(m^2+2mn+m^2n^2)(mn-n^2+2m^2n^2)c+2(m^2+2mn+m^2n^2)^2m^2n^2\ge 0.$$
Recall that
$$m^2+2mn+m^2n^2=(1+mn+\sqrt{1-m^2})(1+mn-\sqrt{1-m^2}),$$
therefore
the inequality is equivalent to
$$2(1+mn+\sqrt{1-m^2})^2m^2n^2-2(1+mn+\sqrt{1-m^2})(mn-n^2+2m^2n^2)mn+(m^2n^2-n^2)m^2n^2\ge 0.$$
Using
$$(1+mn+\sqrt{1-m^2})^2-2(1+mn)(1+mn+\sqrt{1-m^2})+m^2+2mn+m^2n^2=0,$$
we see that the last inequality is equivalent to
$$mn^2(m+n)(1+mn+\sqrt{1-m^2})\ge \frac{1}{2}m^2n^2(2m^2+n^2+4mn+m^2n^2),$$
which further reduces to
$$2(m+n)\sqrt{1-m^2}\ge (1-m^2)(-2m-2n-mn^2).$$
This is true since $2(m+n)\ge 0\ge -2m-2n-mn^2$, and finally we get
$$2(m+n)\sqrt{1-m^2}\ge 2(m+n)(1-m^2)\ge  (1-m^2)(-2m-2n-mn^2).$$
                      
\end{proof}

\begin{rmk}
It might seem natural to define $\Phi_{m,n}$ with $c=-\frac{1}{2}mn(1+mn+\sqrt{1-\max\{m^2,n^2\}})$. However, then $\Phi_{m,n}$ may not be positive. For example, for $m=\frac{25}{32}, n=-\frac{1}{2}$ and $x=(i,1,1)$  we get $\det \Phi_{m,n}(\x\x^{\ast})<0$.

\end{rmk}



\end{document}